\documentclass[11pt,oneside]{article}
\usepackage{amsmath}
\usepackage{verbatim}
\usepackage{graphicx}
\usepackage{amsfonts}
\usepackage{algorithmic}
\usepackage{rotating}
\usepackage{multicol}
\usepackage{enumerate}
\usepackage{pdflscape}
\usepackage{xcolor}
\usepackage[normalem]{ulem}

\usepackage{bm}

\usepackage[ruled, linesnumbered, vlined]{algorithm2e}

\usepackage{longtable}
\usepackage{natbib}
\usepackage{booktabs}
\usepackage[T1]{fontenc}
\usepackage{graphics}
\usepackage{amssymb}
\usepackage{pgf}
\usepackage{colortbl}
\usepackage{multirow}
\usepackage{booktabs}
\usepackage{float,adjustbox,threeparttable}
\usepackage{appendix}

\usepackage{amsthm}

\newtheorem{proposition}{Proposition}

\def\squareforqed{\hbox{\rlap{$\sqcap$}$\sqcup$}\vspace{2mm}}
\def\qed{\ifmmode\else\unskip\quad\fi\squareforqed}
\def\smartqed{\def\qed{\ifmmode\squareforqed\else{\unskip\nobreak\hfil
\penalty50\hskip1em\null\nobreak\hfil\squareforqed
\parfillskip=0pt\finalhyphendemerits=0\endgraf}\fi}}

\smartqed

\newtheorem{theorem}{Theorem}
\newenvironment{myproof}[1][Proof]{\noindent\textbf{#1.} }{\qed}

\newcommand{\DEFINE}{\overset{def.}{\ \equiv\ }}

\newcommand{\EVIDENZIA}[2][black]{\textcolor{#1}{#2}}
\newcommand{\SECREF}[1]{\EVIDENZIA[black]{Section \ref{#1}}}
\newcommand{\APXREF}[1]{\EVIDENZIA[black]{Appendix \ref{#1}}}
\newcommand{\ALGREF}[1]{\EVIDENZIA[black]{Algorithm \ref{#1}}}
\newcommand{\TEOREF}[1]{\EVIDENZIA[black]{Theorem \ref{#1}}}
\newcommand{\PROREF}[1]{Proposition \ref{#1}}

\newcommand{\EQUREF}[1]{\EVIDENZIA[black]{\eqref{#1}}}

\newcommand{\LINEREF}[1]{\EVIDENZIA[black]{\ref{#1}}}
\newcommand{\FACTREF}[1]{\EVIDENZIA[black]{\eqref{#1}}}

\newcommand{\textVT}[1]{\ensuremath{\bm{#1}}}
\newcommand{\XBF}{\textVT{x}}
\newcommand{\YBF}{\textVT{y}}
\newcommand{\ZBF}{\textVT{z}}
\newcommand{\SBF}{\textVT{s}}
\newcommand{\SHBF}{\ensuremath{\hat{\textVT{s}}}}
\newcommand{\lpBF}{\ensuremath{\textVT{l}^{+}}}
\newcommand{\lmBF}{\ensuremath{\textVT{l}^{-}}}
\newcommand{\DBF}{\textVT{d}}
\newcommand{\qBF}{\textVT{q}}
\newcommand{\QBF}{\textVT{Q}}
\newcommand{\OBF}{\ensuremath{\boldsymbol{0}}}
\newcommand{\XBFI}[1][n]{\ensuremath{x^{[#1]}}}
\newcommand{\DBFI}[1][n]{\ensuremath{\delta^{[#1]}}}

\newcommand{\STATION}{station}

\newcommand{\BOXED}[1]{\framebox[1.1\width]{\ensuremath{#1}}}

\newcommand{\NUMV}{\ensuremath{w}}

\newcommand{\ENDEPOCH}{\ensuremath{\overline{e}}}
\newcommand{\STARTEPOCH}{\ensuremath{\underline{e}}}
\newcommand{\SMINT}[1]{\ensuremath{\underline{s}_{[\STARTEPOCH,#1]}}}
\newcommand{\SMAXT}[1]{\ensuremath{\overline{s}_{[\STARTEPOCH,#1]}}}

\newcommand{\UP}{\ensuremath{\beta}}
\newcommand{\DOWN}{\ensuremath{\alpha}}

\newcommand{\ENDSTOCK}[1]{\ensuremath{s^{(\ENDEPOCH)}(#1)}}

\newcommand{\LI}[2]{\ensuremath{L(#2\ |\ #1)}}  

\newcommand{\LQSTAR}{\ensuremath{L^{Q}}}
\newcommand{\LISTAR}{\ensuremath{L^{\infty}}}

\newcommand{\XQSTAR}{\ensuremath{x^{Q}}}
\newcommand{\XQUSTAR}{\ensuremath{\overline{x}^{Q}}}
\newcommand{\XQDSTAR}{\ensuremath{\underline{x}^{Q}}}
\newcommand{\XQRSTAR}{\ensuremath{[\XQDSTAR,\XQUSTAR]}}
\newcommand{\XQSSTAR}{\ensuremath{\bm{X}^{Q}}}

\newcommand{\XISTAR}{\ensuremath{x^{\infty}}}
\newcommand{\XIUSTAR}{\ensuremath{\overline{x}^{\infty}}}
\newcommand{\XIDSTAR}{\ensuremath{\underline{x}^{\infty}}}
\newcommand{\XIRSTAR}{\ensuremath{[\XIDSTAR,\XIUSTAR]}}
\newcommand{\XISSTAR}{\ensuremath{\bm{X}^{\infty}}}

\newcommand{\XVECTrow}[3][x]{\ensuremath{[#1_{#2},\ldots,#1_{#3}]}}

\newcommand{\INISTOCKEPOCHBYVECT}[2][x]{\ensuremath{s^{(e_{#2}-1)}(#1)}}

\newcommand{\SX}[2][x]{\ensuremath{s^{(e_{#2}-1)}(\bm{#1})}}
\newcommand{\SXPIU}[2][x]{\ensuremath{s^{(e_{#2+1}-1)}(\bm{#1})}}

\newcommand{\SZERO}{\ensuremath{s^{(0)}}}

\newcommand{\PROBLEMA}[1]{\ensuremath{\mathcal{P}(#1)}}
\newcommand{\PROBLEMAAUG}[2][\delta]{\ensuremath{\bar{\mathcal{P}}(#2,#1)}}

\newcommand{\LVECTBKW}[2][x]{\ensuremath{L_{J_{#2}}\left(\XVECTrow[#1]{1}{#2}\right)}}
\newcommand{\LVECTBKWONE}[1][x]{\ensuremath{L_{J_{1}}\left([#1_1]\right)}}
\newcommand{\LVECTBKWO}[1]{\LVECTBKW[0]{#1}}
\newcommand{\PIOTTBKW}[1]{\ensuremath{\min_{\XVECTrow[x]{1}{#1}}\LVECTBKW{{#1}}}}

\newcommand{\AUGMINTI}[1][i]{\ensuremath{\bar{I}_{#1}}}
\newcommand{\AUGMINTJ}[1][i]{\ensuremath{\bar{J}_{#1}}}

\newcommand{\LVECTBKWXAUGMN}[2]{\ensuremath{L_{[\AUGMINTJ[#2],\delta]}\left(\XVECTrow[#1]{1}{#2}\right)}}
\newcommand{\LVECTBKWAUGMN}[1]{\LVECTBKWXAUGMN{x}{#1}}
\newcommand{\LVECTBKWOAUGMN}[1]{\LVECTBKWXAUGMN{0}{#1}}

\newcommand{\PIOTTBKWAUGMN}[1]{\min_{\XVECTrow[x]{1}{#1}}\LVECTBKWAUGMN{{#1}}}

\newcommand{\LIXS}[3][x]{L_{I_{#2}}\left(#1\ \left|\ #3\right.\right)}

\newcommand{\LIOTTS}[3][x]{\ensuremath{L_{I_{#2}}\left(#1_{#2}\ \left|\ #3\right.\right)}}  

\newcommand{\LIOTT}[2][x]{\LIOTTS[#1]{#2}{\SX[#1]{{#2}}}}
\newcommand{\LIOTTGIVENS}[3][x]{\LIOTTS[#1]{#2}{#3}}

\newcommand{\LIOTTAUGMN}[3][x]{L_{[\AUGMINTI[#2],\delta]}\left(#1_{#2}\ \left|\ #3\right.\right)}
\newcommand{\LIOTTAUGMNX}[3][x]{L_{[\AUGMINTI[#2],\delta]}\left(#1\ \left|\ #3\right.\right)}

\newcommand{\HZEROX}{\underline{h}}
\newcommand{\SMINBEFOREHZEROX}{\underline{s}}

\newcommand{\ONEINTPROBLEM}{\textit{1-Intervention}}

\newcommand{\LSTAR}[1][i]{\ensuremath{L^{*}_{J_{#1}}}}
\newcommand{\LSTARAUGMN}[1][i]{\ensuremath{L^{*}_{[\AUGMINTJ[#1],\delta]}}}
\newcommand{\XSTARAUGMN}{\ensuremath{\bar{\XBF}}}

\newcommand{\AlgoVehicleIntervention}{Vehicle-Intervention}
\newcommand{\CiteAlgoVehicleIntervention}{\normalfont{\texttt{\AlgoVehicleIntervention}}}
\newcommand{\AlgoGlobalBackward}{Global-Backward}
\newcommand{\CiteAlgoGlobalBackward}{\normalfont{\texttt{\AlgoGlobalBackward}}}
\newcommand{\AlgoGlobalBackwardRec}{Global-Backward-Recursive}
\newcommand{\CiteAlgoGlobalBackwardRec}{\normalfont{\texttt{\AlgoGlobalBackwardRec}}}

\begin{document}

\title{The one-station bike repositioning problem}
\author{\it E. Angelelli$^1$, A. Mor$^2$, M.G. Speranza$^1$ \\
{\small \it $^1$Department of Economics and Management}\\
{\small \it University of Brescia, Italy}\\
{\small \{enrico.angelelli, grazia.speranza\}@unibs.it}\\
{\small \it $^2$Department of Management, Economics and Industrial Engineering}\\
{\small \it Politecnico di Milano, Italy}\\
{\small andrea.mor@polimi.it}\\
\textbf{Accepted for publication in Discrete Applied Mathematics}
}

\date{}
\maketitle

\begin{abstract}

In bike sharing systems the quality of the service to the users strongly depends on the strategy adopted to reposition the bikes. The bike repositioning problem is in general very complex as it involves different interrelated decisions: the routing of the repositioning vehicles, the scheduling of their visits to the stations, the number of bikes to load or unload  for each station and for each vehicle that visits the station.
In this paper we study the problem of optimally loading/unloading vehicles that visit the same station at given time instants of a finite time horizon.
The goal is to minimize the total lost demand of bikes and free stands in the station.
We model the problem as a mixed integer linear programming problem and present an optimal algorithm that runs in linear time in the size of the time horizon.

\vspace{3mm}

\end{abstract}

\section{Introduction}

In a Bike Sharing System (BSS) one of the most important decisions at the operational level is the adoption of a bike repositioning strategy to mitigate the imbalances caused by user demand over time and space. These strategies can be divided in user-based and vehicle-based. In the former, reward systems are put in place to encourage a use of the system that counteracts the imbalance. In the latter, one or more vehicles are deployed to move bike to and from stations to rebalance the bike inventory. In this paper we consider the latter and focus on the operations of a fleet of vehicles in a single station.

Vehicle-based bike repositioning strategies can be further classified as static or dynamic.
When a static bike repositioning strategy is adopted, the bike repositioning operations are performed when the user demand is low
compared to the rest of the day, typically at night. In this case, real-time user demand is disregarded.
Dynamic bike repositioning is performed during the day, accounting for the real-time status of the system and the behavior of the users, as well as its forecast.
Typically, the dynamic problem is tackled by solving a static time-dependent subproblem at the beginning of the operating day based on the current state and a forecast of the demand during the day.
Decisions are then periodically re-evaluated at later times  with updated information.

Various objectives have been considered in the literature to evaluate the quality of the relocation solution. The most common one is the
maximization of the quality of service, typically measured as the number of lost withdraws (due to the lack of bikes) and/or returns (due to the lack of stands). The former causes discomfort as the user has to either wait for a bike to be returned at the desired departure station, move to another station where bikes are available, or resort to another means of transportation. In the latter case, the user has to cycle to a different station with free stands or, if available, resort to alternative ways to return the bike (e.g., with a lock issued by the BSS operator).

Decisions to be taken in a repositioning strategy include those about
the amount of bikes to be loaded or unloaded in each station,
the routing of the vehicles,
and, in the dynamic setting, the timing of the visits of the vehicles at each station.
We collect the models that approach these decisions under the name of bike repositioning problems (BRPs).
These decisions are interrelated.
The routing determines the sequence of the stations visited by the vehicles and therefore constrains their arrival time at each station. The number of bikes present at each station depends on the arrival time.

The decisions taken for the operations of a vehicle in a station influence (interact with) not only the actions of the same vehicle in the following stations along its route but also the actions of other vehicles that could be scheduled to visit the same station.
In fact, the interaction between the operations of multiple vehicles is one of the main challenges when solving a problem in the class of BRPs.
The need of multiple visits, by the same vehicle or different vehicles, to the same station is a consequence of the limited capacity of the repositioning vehicles, and of the dynamic nature of the user demand. Both loading and unloading operations need to be performed across the stations and across time at the same station.

When multiple visits to a station are considered, the interaction among the vehicles should be taken into account when planning the loading/unloading operations of each vehicle.
Any decision taken for a vehicle has an impact on the operations that the following vehicles can perform along their routes.
For example, if a vehicle loads bikes to tackle a shortage of stands, and the next vehicle unloads bikes to prevent a later shortage of bikes,
the number of bikes loaded by the first vehicle may have an impact on the size of the subsequent  bike shortage.
In other words, loading an exceedingly large number of bikes on the first vehicle may cause a worsening of the shortage of bikes at a later time that the second vehicle may be unable to cope with.

\paragraph{Contribution.} Drawing from the experience gained in a collaboration with Brescia Mobilità SpA, operating the station-based BSSs in Brescia, Italy (see \cite{angelelli2022simulation}), in this paper we study the problem of optimally loading/unloading vehicles that visit the same station at given time instants of a discrete and finite time horizon.
The net users demand at the station is known (provided by the forecast) at each time instant of the time horizon.
The goal is to minimize the total lost demand of bikes and free stands.
The station is characterized by an initial inventory stock and a known capacity.
Likewise, each vehicle is characterized by a capacity and an initial load of bikes.
We model the optimization problem, that accounts for the interaction of the operations of the vehicles, as a mixed integer linear programming problem whose linear relaxation is shown to always have an integer solution. We then present a solution algorithm that runs in linear time in the size of the time horizon.
Given the forecast of net users demand on each instant (epoch) of the time horizon, the algorithm finds the type of operation (loading/unloading) and the number of bikes (to be loaded/unloaded) for each vehicle that visits the station to minimize the number of unsatisfied requests, taking into account the limits imposed by the capacity and the load at the time of the visit.
It also finds the number of requests that remain unsatisfied  regardless of the load and capacity of the vehicles. Therefore, the algorithm could be also used to calibrate the load of the vehicles upon their arrival at a station and to evaluate the effectiveness of certain visiting times at a station.
The problem addressed can be seen as a subproblem to be solved in a solution approach to a more general repositioning problem. The availability of a very efficient optimal algorithm for  the solution of the subproblem will contribute to the efficiency of the approach.
In this sense our work  provides a contribution to the design of solution approaches to more general models in the class of BRPs.

\paragraph{Literature Review.}

The design and management of a BSS involve many issues that have been tackled by the optimization literature (see \cite{laporte2018shared,shui2020review,vallez2021challenges}). Among these problems, those regarding the mitigation of spatio-temporal imbalances in the demand and offer of bikes by the users have received a considerable attention.
Many works consider a static setting, where bike relocation is assumed to take place at night, when user demand is negligible (see \cite{benchimol2011balancing,chemla2013bike,dell2014bike,erdougan2014static}). In this context, the decisions can be decomposed in those regarding the inventory of the station and those about the routing of the relocating vehicles.
Considering the nighttime-only relocation of bikes can be a limitation in systems experiencing high fluctuations in user demand. This limitation is overcome in the dynamic setting where the relocation takes place during the day and therefore decisions can be revised multiple times. In this setting, the decisions to be taken include the inventory of each station, and the scheduling and routing of the vehicles performing the relocation must be taken. These decisions have been tackled both separately (e.g., \cite{regue2014proactive}) and in an integrated way (e.g., \cite{zhang2017time,ghosh2017dynamic,zheng2021repositioning}).
It must be noted that the dynamic setting requires to accurately forecast user behavior. This requirement is considered in various contributions (e.g., \cite{regue2014proactive,alvarez2016optimizing,yoshida2019practical,lee2020optimal}).

As noted above, in this paper we aim at tackling the problem of coordinating the operations of multiple vehicles visiting a single station over time.
A summary of the literature for this class of BRPs that is relevant for our contribution is presented in Table \ref{tab:litrev}. For each reference, the table classifies the papers according to: the time setting of the problem (static or dynamic), the number of vehicles for the repositioning of bikes, whether stations are allowed to be visited multiple times, and whether vehicle interaction is considered or not.
Considering the real-time evolution of the system and the forecast of users behavior makes the bike repositioning a very difficult problem to solve. It is therefore not surprising that early contributions generally considered the static setting. Furthermore, a considerable portion of the literature examined considers the case of a single vehicle or does not allow for multiple visits or vehicle interaction.
To further highlight the difficulty of the task, it is also worth noting that, in many of the papers, the fact that stations might be visited multiple times is not the result of a repositioning plan that explicitly considers this possibility and the consequent vehicle interactions.
In these papers this is the result of either the problem being modeled as a Markov Decision Process where the set of actions includes the visit to one among all the stations of the system, or the problem being solved by considering a rolling horizon scheme where the reoptimization step considers a visit to all the stations, that is, without removing from the set of candidate stops stations that have been visited at earlier times.

We report that, to the best of our knowledge, the only work considering the problem on a single station is \cite{raviv2013optimal}. The authors model the user requests to rent or return a bike as stochastic processes and introduce a dissatisfaction function, establishing its convexity and devising an accurate and efficient approximation for its estimation.

\bigskip
The paper is organized as follows. First, in \SECREF{sec:probDef} the problem is defined. Then, in \SECREF{sec:localOpt} the optimization problem is discussed for specific subsets of the time horizon. The results are then exploited in \SECREF{sec:globalOpt} to present a solution algorithm for the problem on the entire time horizon.

\begin{table}[H]
    \centering
    \begin{adjustbox}{max size={1\textwidth}{0.38\textheight}}
    \begin{threeparttable}
        \begin{tabular}{lllcc}
        \hline
            \textbf{Reference} & \begin{tabular}[r]{@{}c@{}} \textbf{Time} \\ \textbf{setting} \end{tabular} & \textbf{\# of vehicles} & \begin{tabular}[r]{@{}c@{}} \textbf{Multiple} \\ \textbf{visits} \end{tabular} & \begin{tabular}[r]{@{}c@{}} \textbf{Vehicles} \\ \textbf{interaction} \end{tabular} \\ \hline
            \cite{benchimol2011balancing} & Static & Single & - & - \\
            \cite{chemla2013bike} & Static & Single & - & - \\
            \cite{erdougan2014static} & Static & Single & - & - \\
            \cite{ho2014solving} & Static & Single & - & - \\
            \cite{li2016multiple} & Static & Single & - & - \\
            \cite{szeto2016chemical} & Static & Single & - & - \\
            \cite{cruz2017heuristic} & Static & Single & \checkmark & - \\
            \hline

            \cite{lin2012geo} & Static & Multiple & - & - \\
            \cite{gaspero2013hybrid} & Static & Multiple & - & - \\
            \cite{forma20153} & Static & Multiple & - & - \\
            \cite{dell2016destroy} & Static & Multiple & - & - \\
            \cite{szeto2018exact} & Static & Multiple & - & - \\
            \cite{schuijbroek2017inventory} & Static & Multiple\tnote{1} & - & - \\
            \cite{bulhoes2018static} & Static & Multiple & \checkmark & - \\
            \cite{rainer2013balancing} & Static & Multiple & \checkmark & \checkmark \\
            \cite{raviv2013static} & Static & Multiple & \checkmark & \checkmark \\
            \cite{angeloudis2014strategic} & Static & Multiple & \checkmark & \checkmark \\
            \cite{alvarez2016optimizing} & Static & Multiple & \checkmark & \checkmark \\
            \cite{ho2017hybrid} & Static & Multiple & \checkmark & \checkmark \\
            \cite{wang2018static} & Static & Multiple & \checkmark & \checkmark \\
            \cite{casazza2018multiple} & Static & Multiple & \checkmark & \checkmark\tnote{4} \\
            \hline
            \cite{caggiani2012modular} & Dynamic & Single & - & - \\
            \cite{shui2018dynamic} & Dynamic & Single & \checkmark & - \\
            \cite{brinkmann2019dynamic} & Dynamic & Single & \checkmark\tnote{2} & - \\
            \hline
            \cite{kloimullner2014balancing} & Dynamic & Multiple & - & - \\
            \cite{huang2022monte} & Dynamic & Multiple\tnote{1} & - & - \\
            \cite{legros2019dynamic} & Dynamic & Multiple\tnote{1} & \checkmark\tnote{2} & - \\
            \cite{brinkmann2020multi} & Dynamic & Multiple & \checkmark\tnote{2} & - \\
            \cite{contardo2012balancing} & Dynamic & Multiple & \checkmark\tnote{3} & - \\
            \cite{zhang2017time} & Dynamic & Multiple & \checkmark\tnote{3} & - \\
            \cite{chiariotti2020bike} & Dynamic & Multiple & \checkmark\tnote{3} & - \\
            \cite{angelelli2022simulation} & Dynamic & Multiple & \checkmark\tnote{3} & - \\
            \cite{ghosh2017dynamic} & Dynamic & Multiple & \checkmark & \checkmark \\
            \cite{lowalekar2017online} & Dynamic & Multiple & \checkmark &  \checkmark \\
            \cite{ghosh2019improving} & Dynamic & Multiple & \checkmark & \checkmark \\
            \cite{zheng2021repositioning} & Dynamic & Multiple & \checkmark & \checkmark \\ \hline
    \end{tabular}
    \begin{tablenotes}
    \item[1] Stations are clustered and each cluster is assigned to one vehicle.
    \item[2] A Markov Decision Process is presented where routing decisions involve only one vehicle visit per station. Further visits to a station are allowed only in subsequent stages.
    \item[3] Each station can be served at most once by one vehicle within a repositioning time-window. A rolling horizon scheme is presented, with multiple repositioning time-windows, meaning that a station can be served in multiple time-windows during the day.
    \item[4] The vertices can be visited several times and by distinct vehicles, but temporary storage is not allowed. This means that when inventory units are delivered to some vertices, they cannot be picked up later again. 
    \end{tablenotes}
    \end{threeparttable}
    \end{adjustbox}
    \caption{Literature contributions with respect to the characteristics relevant for this work.}
    \label{tab:litrev}
    \end{table}

\section{Notation and problem definition} \label{sec:probDef}

Let us consider a bike-station (\textit{\STATION} from now on) with known capacity $C$ and initial stock level $\SZERO\in[0,C]$ at time $t^{(0)}$.
Let the time horizon consist of $m$ time instants $t^{(h)}$ indexed with $h\in T=\{1,\ldots,m\}$, where $t^{(h-1)}<t^{(h)}$ for all $h=1,\ldots,m$, and let us call \textit{epochs} the indices $h\in T$.
The \STATION\ receives at epoch $h$ a known request of \textit{net flow} of bikes (difference between rentals and returns) $d^{(h)}$, where $d^{(h)}>0$ represents an incoming net flow and $d^{(h)}<0$ represents an outgoing net flow. 

A set $V=\{v_1,\ldots,v_\NUMV\}$ of $\NUMV$ vehicles is available to visit the \STATION\ with a fixed schedule.
More precisely, vehicle $v_i$ visits the \STATION\ at a fixed epoch $e_i\in T$. We denote by $H=\{e_1,\ldots,e_\NUMV\}\subseteq T$ the set of visiting epochs and w.l.o.g. we assume that  $e_1=1$ and vehicles are indexed according to visiting epochs so that $e_i<e_{i+1}$ for $i=1,\ldots,\NUMV-1$. Every vehicle $v_i$ has a fixed capacity $Q_i$ and a known load $q_i\in[0,Q_i]$ at visit epoch $e_i$.

Thus, at epoch $e_i$ the operator of vehicle $v_i$ may move, in either direction, some bikes between the \STATION\ and the vehicle.
The amount moved at epoch $e_i$ is indicated with $x_i$, where $x_i>0$ represents a bike flow from the vehicle to the \STATION\ (i.e., bike unloading from the vehicle) and $x_i<0$ represents a bike flow from the \STATION\ to the vehicle (i.e., bike loading on the vehicle). We call this amount an \textit{intervention}.
Intervention $x_i$ is bounded from above by the load $q_i$ of the vehicle and from below by its residual capacity $q_i-Q_i$.
Thus, the intervention is feasible if and only if
\begin{align}\label{eq01:update_q}
    x_i\in[q_i-Q_i,q_i]\quad \text{(or, equivalently $q_i-x_i\in[0,Q_i]$).}
\end{align}

As the intervention $x_i$ and the request $d^{(e_i)}$ take place at epoch $e_i$ the final result is a (post intervention) virtual stock given by
\begin{align}\label{eq02:update_shat_int}
    \hat{s}^{(e_i)}=s^{(e_i-1)}+d^{(e_i)}+x_i.
\end{align}
In any other epoch $h\notin H$ the virtual stock is given by
\begin{align}\label{eq03:update_shat_noint} 
    \hat{s}^{(h)}=s^{(h-1)}+d^{(h)}.
\end{align}
At any epoch $h\in T$, if the virtual stock is greater than the \STATION\ capacity, or below zero, the difference is accounted for as a loss.
More precisely, we define the \textit{surplus loss} as
\begin{align}\label{eq04:update_lsurplus}  
    l^{(h)+}=\max(0,\hat{s}^{(h)}-C)
\end{align}
and the \textit{stockout loss} as
\begin{align}\label{eq05:update_lstockout}   
    l^{(h)-}=\max(0,-\hat{s}^{(h)}).
\end{align}
It is clear that at most one of these two quantities may be positive.
The actual (post intervention) stock level in the \STATION\ is
\begin{align}\label{eq06:update_s}
    s^{(h)}=\hat{s}^{(h)}-l^{(h)+}+l^{(h)-}
\end{align}
where by construction $s^{(h)} \in[0,C]$, $s^{(h)}=0$ if $l^{(h)-}>0$ and $s^{(h)}=C$ if $l^{(h)+}>0$. We define the total amount of lost requests as $L=\sum_{h=1}^{m}\left(l^{(h)+}+l^{(h)-}\right)$.

\begin{figure}[H]
	\centering
	\includegraphics[width=0.65\linewidth]{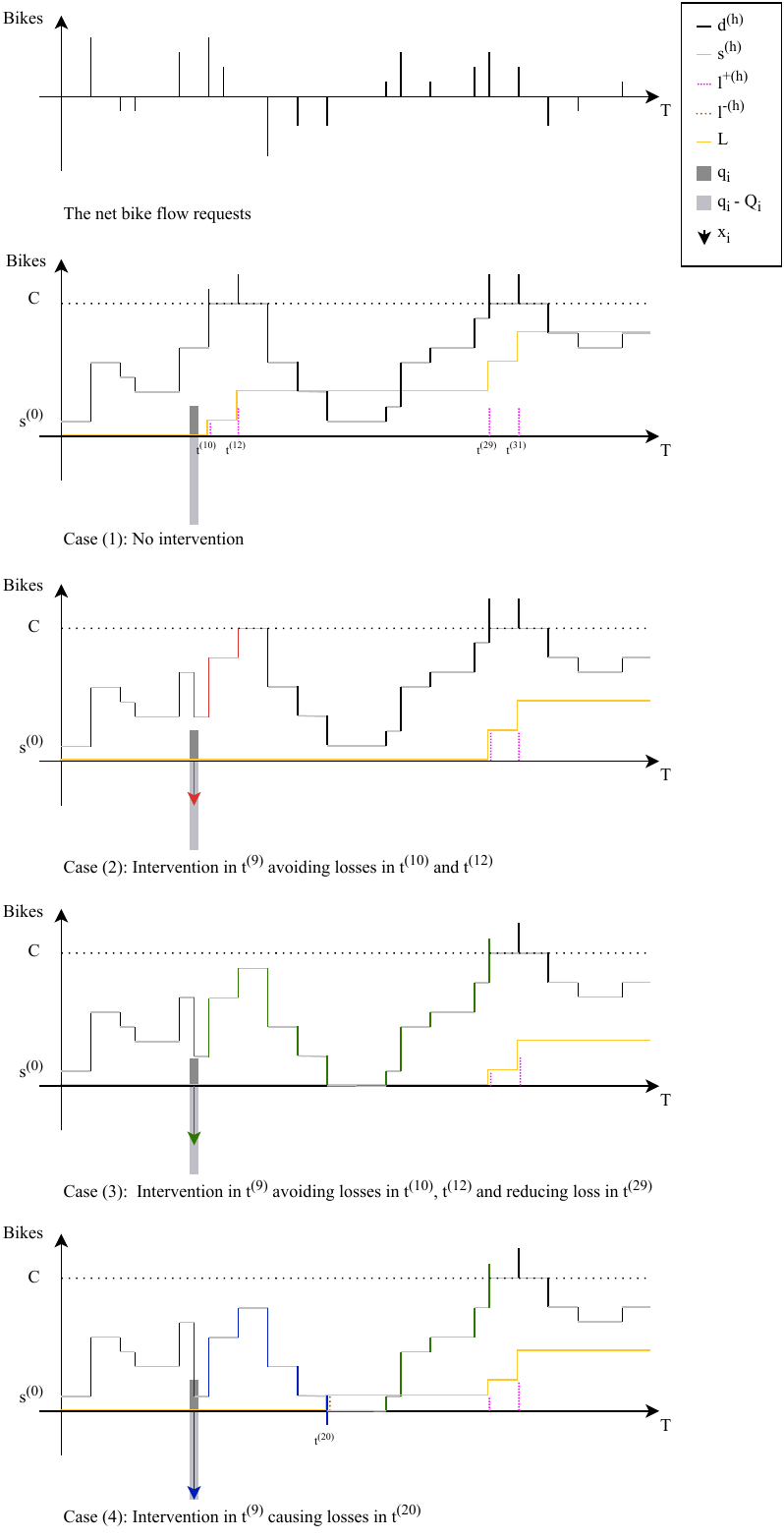}
	\caption{A representation of the quantities involved in the problem: the net flow of bikes (above), and resulting actual stock, and surplus and stockout losses. The effect on $s^{(h)}$, $l^{(h)-}$, and $l^{(h)-}$ of different values for an intervention in $t^{(9)}$.} 
	\label{fig:fig1}
\end{figure}

Figure \ref{fig:fig1} reports an example of the net flow of bikes $\DBF$. A vehicle $v_9$ with capacity $Q_9$ and load $q_9$ is visiting the station at time $t^{(9)}$. Note that, in the figure, a small jitter is present to avoid superimposition. For a given initial stock level $\SZERO$, Case (1) represents the resulting stock $s^{(h)}$, surplus loss $l^{(h)+}$, and stockout loss $l^{(h)-}$ (for $h \in T$) in the case of null intervention. The total loss $L$ is also reported. While no temporal characterization is given for $L$, its value over time is represented in Figure \ref{fig:fig1} for ease of visualization.
It is worth pointing out that losses might occur in an epoch $h$ when an intervention is not allowed ($h\notin H$), unless they are avoided by an earlier intervention.
Observe that earlier interventions should take into account that the \STATION\ stock level is waving up and down due to the net flow, and since an intervention at epoch $h$ produces a shift in the stock level from epoch $h$ onward, it may avoid some losses in the next epochs, but may also generate new losses in one or more of the next epochs.
This is exemplified in Case (2), where the intervention planned for $t^{(9)}$ avoids the losses occurring in $t^{(10)}$ and $t^{(12)}$ when no interventions can be performed. Increasing the number of bikes loaded in $t^{(9)}$ reduces the loss occurring in $t^{(29)}$ (Case (3)), but because of the waving of the net demand, a larger intervention induces a stockout loss in $t^{(20)}$ (Case (4)), without modifying the actual stock at later times (i.e., the loss at $t^{(29)}$ and later is unchanged w.r.t. Case (3)).

Problem instance is defined by parameters \SZERO, $C$ and vectors $\DBF=[d^{(1)},\ldots,d^{(m)}]$, $\qBF=[q_1,\ldots,q_\NUMV]$, $\QBF=[Q_1,\ldots,Q_\NUMV]$. Decision variables are described by a vector of interventions $\XBF= \XVECTrow{1}{\NUMV}$. 
When required by the discussion, we will refer to quantities $s^{(h)}$, $\hat{s}^{(h)}$, $l^{(h)+}$, and $l^{(h)-}$ as $s^{(h)}(\XBF)$, $\hat{s}^{(h)}(\XBF)$, $l^{(h)+}(\XBF)$, and $l^{(h)-}(\XBF)$, respectively, to explicitly express their dependence from interventions $\XBF$. Accordingly, the total amount of lost requests $L$ will be denoted by $L(\XBF)=\sum_{h=1}^{m}\left(l^{(h)+}(\XBF)+l^{(h)-}(\XBF)\right)$.

Thus, the problem can be formulated as a MILP as follows
\begin{align}
\min \quad & L(\XBF) &\label{eq07:base:L} \\
\nonumber s.t. \quad & \text{equations \EQUREF{eq01:update_q}-\EQUREF{eq06:update_s}} & \text{for } h=1,\ldots,m.\\
& x_i \in \mathbb{Z} & \text{for } i=1,\ldots,\NUMV.\label{eq08:base:Z}
\end{align}

The linear relaxation of this problem is obtained by removing the integrality conditions $x_i \in \mathbb{Z}$. It can be shown that the problem always has an integer solution.

\begin{theorem}
An integer solution always exists to the linear relaxation of Problem \eqref{eq01:update_q}-\eqref{eq08:base:Z}.

\end{theorem}
\begin{myproof}
The claim follows from observing that the matrix of the coefficients of the variables in the constraints is totally unimodular. This is true because
\begin{itemize}
\item every entry is $0$, $+1$, or $-1$;
\item every column contains at most two non-zero (i.e., $+1$ or $-1$) entries;
\item each variable that appears in more than one constraint, that is each $s^{(h)}$ and each $\hat{s}^{(h)}$, with $h\in T$, appears with coefficients $+1$ and $-1$.
\end{itemize}
These are sufficient conditions for the matrix to be totally unimodular (see \cite{schrijver2003combinatorial}).
\end{myproof}

In this paper we show that the problem can be solved in $O(m)$ time.

Observe that $\XBF=\OBF$ is a feasible vector of interventions provided $q_i\in[0,Q_i]$ for all $v_i\in V$. Observe also that the total loss is minimized by maximizing the quantity $L(\OBF)-L(\XBF)$ which we call recovered (or saved) loss.
More in detail, if $(l^{(h)+}(\OBF)+l^{(h)-}(\OBF))>(l^{(h)+}(\XBF)+l^{(h)-}(\XBF))$ we say that interventions \XBF\ recovers a loss at epoch $h$; while, if $(l^{(h)+}(\OBF)+l^{(h)-}(\OBF))<(l^{(h)+}(\XBF)+l^{(h)-}(\XBF))$ we say that interventions \XBF\ creates (or increments) a loss at epoch $h$.

Note that, for the sake of ease of reading, all the proofs of the propositions are provided in \APXREF{app:proof}. For the same reason, the appendix also reports a table of notation.

\section{Interval optimization} \label{sec:localOpt}

In this section we study the case of the problem of a single intervention to be taken at epoch $\STARTEPOCH$ (i.e., $H=\{\STARTEPOCH\}$) by a vehicle with capacity $Q$ and load $q$ to minimize the loss realized within an epoch interval $[\STARTEPOCH,\ENDEPOCH]$. 
We call this particular problem the \ONEINTPROBLEM\ problem. As we have a single decision, the vector of interventions \XBF\ reduces to a scalar value $x$.

We focus first on the uncapacitated version of the \ONEINTPROBLEM\ problem where constraint \EQUREF{eq01:update_q} is relaxed and characterize the set of optimal interventions. Then, we establish the relationship with the capacitated problem and design a linear time optimal algorithm.
To characterize optimal interventions, the sufficient conditions to improve an intervention are described (\SECREF{sec:ImprovingInterventions}). Then, sufficient conditions for the worsening of an intervention are defined (\SECREF{ref:WorseningInterventions}). These results allow for the definition of the conditions for the optimality of an intervention (\SECREF{sec:OptimalityConditions}) and, eventually, the definition of an algorithm running in $O(\ENDEPOCH-\STARTEPOCH)$ time (\SECREF{sec:intervaloptimization}).

Some additional notation is now introduced for this particular \ONEINTPROBLEM\ version of the problem.
Given an intervention $x$ at epoch $\STARTEPOCH$ we define

\begin{itemize}
    \item $\SMINT{h}(x)=\min(s^{(k)}(x)\ |\ \STARTEPOCH\leq k\leq h)$, that is the minimum stock level in epoch interval $[\STARTEPOCH,h]$;
    \item $\SMAXT{h}(x)=\max(s^{(k)}(x)\ |\ \STARTEPOCH\leq k\leq h)$, that is the maximum stock level in epoch interval $[\STARTEPOCH,h]$;
    \item $h^{(+)}(x) = \min(k\in[\STARTEPOCH,\ENDEPOCH]|\ l^{(k)+}(x)>0)$ if there is an epoch with surplus loss, $+\infty$ otherwise;
    \item $h^{(-)}(x) = \min(k\in[\STARTEPOCH,\ENDEPOCH]\ |\ l^{(k)-}(x)>0)$ if there is an epoch with stockout loss, $+\infty$ otherwise.
\end{itemize}

\subsection{Improving interventions} \label{sec:ImprovingInterventions}

Here we state sufficient conditions to improve a given intervention $x$.
Considering the case that, for a given intervention $x$, the first loss is of type surplus (i.e., $h^{(+)}(x)<h^{(-)}(x)$), we would like to move from the \STATION\ to the vehicle a quantity $\delta$ so that the stock in the interval $[\STARTEPOCH,h^{(+)}(x)-1]$ is downshifted and the surplus at epoch $h^{(+)}(x)$ is reduced or even zeroed without changing the stock since then.
However, the downshift should not be too large in order to avoid the creation of stockout losses that would also stop the surplus recovery.
In case we manage to zero the loss at epoch $h^{(+)}(x)$, the epoch of first loss is moved forward and we can look for a further improvement of the new intervention (see \SECREF{sec:intervaloptimization}).
A similar argument can be used when $h^{(-)}(x)<h^{(+)}(x)$.

\begin{proposition}[Downshift improvement]\label{PN1:ImproveByDownShift}
Let $x$ be an intervention such that $\overline{h}=h^{(+)}(x)<\infty$ (at least one surplus loss) and define
$\delta^*=\min\left(l^{(\overline{h})+}(x),\ \SMINT{\overline{h}-1}(x)\right)$.

\medskip\noindent
Then, for any $\delta\in[0,\delta^*]$, following properties hold for intervention $x'=x-\delta$
\begin{align*}
    \begin{array}{rlrl}
    l^{(k)-}(x')=l^{(k)-}(x)  & \text{for } k\in [\STARTEPOCH,\ENDEPOCH]
    \\
    l^{(k)+}(x')=l^{(k)+}(x)-\delta  & \text{for } k=\overline{h} & l^{(k)+}(x')=l^{(k)+}(x) & \text{for } k\neq\overline{h}
    \\
    s^{(k)}(x')=s^{(k)}(x)-\delta  & \text{for } k<\overline{h} &  s^{(k)}(x')=s^{(k)}(x) & \text{for } k\geq\overline{h}
    \end{array}
\end{align*}
Moreover, $L(x')=L(x)-\delta$; we also have that $\delta=l^{(\overline{h})+}(x)$ if and only if $h^{(+)}(x) < h^{(+)}(x')$.
\end{proposition}

\medskip
Observe that \PROREF{PN1:ImproveByDownShift} does not provide any improvement in case $\delta^*=0$. 
Otherwise, if $\delta^*>0$, any intervention $x'=x-\delta$ with $\delta\in(0,\delta^*]$ is a strict improvement of $x$, indeed the surplus loss at epoch $h^{(+)}(x)$ is decreased exactly by $\delta$.
Moreover, the stock level before epoch $h^{(+)}(x)$ is uniformly downshifted by $\delta^*$ while the stock level remains unchanged starting from epoch $h^{(+)}(x)$ (in particular we still have $C$ at epoch $h^{(+)}(x)$).
Finally, the first surplus loss induced by the new intervention is pushed forward if and only if $x'$ recovers the whole loss at epoch $h^{(+)}(x)$.

\begin{proposition}[Upshift improvement]\label{PN2:ImproveByUpShift}
Mirror of \PROREF{PN1:ImproveByDownShift}.

Let $x$ be an intervention such that $\underline{h}=h^{(-)}(x)<\infty$ (at least one stockout loss) and define
$\delta^*=\min\left(l^{(\underline{h})-}(x),\ \SMAXT{\underline{h}-1}(x)\right)$.

\medskip\noindent
Then, for any $\delta\in[0,\delta^*]$, following properties hold for intervention $x'=x+\delta$

\begin{align*}
    \begin{array}{rlrl}
    l^{(k)+}(x')=l^{(k)+}(x)  & \text{for } k\in [\STARTEPOCH,\ENDEPOCH] \\
    l^{(k)-}(x')=l^{(k)-}(x)-\delta  & \text{for } k=\underline{h} & l^{(k)-}(x')=l^{(k)-}(x) & \text{for } k\neq\underline{h}  \\
    s^{(k)}(x')=s^{(k)}(x)+\delta  & \text{for } k<\underline{h} &  s^{(k)}(x')=s^{(k)}(x) & \text{for } k\geq\underline{h}
    \end{array}
\end{align*}
Moreover, $L(x')=L(x)-\delta$; we also have that $\delta=l^{(\underline{h})-}(x)$ if and only if $h^{(-)}(x) < h^{(-)}(x')$.
\end{proposition}

\subsection{Worsening interventions} \label{ref:WorseningInterventions}
Here we state sufficient conditions to say that an intervention $x$ cannot be improved by a downshift $x'<x$ or by an upshift $x'>x$.
Let us consider the case that for a given intervention $x$ the stock is zero before any surplus loss. In this case no downshift could recover the surplus loss as it would generate stockout losses at one or more earlier epochs preventing the downshift at the later epoch.

\begin{proposition}[Worsening by downshift]\label{PN3:WorsenByDownShift}
Let $x$ be an intervention and $\hat{h}=\min(h^{(+)}(x),\ENDEPOCH)$.
If $\SMINT{\hat{h}}(x)=0$ (the stock level goes to zero in the interval and before the first epoch with a surplus loss, if any), then, for all $\delta>0$ the following properties hold for intervention $x'= x-\delta$
\begin{align*}
    & L(x')=L(x)+\delta\\
    & \ENDSTOCK{x'}=\ENDSTOCK{x}
\end{align*}
where the new loss is of type stockout and is realized at one or more epochs earlier than $\hat{h}$.
\end{proposition}

\begin{proposition}[Worsening by upshift]\label{PN4:WorsenByUpShift}
Mirror of \PROREF{PN3:WorsenByDownShift}.
Let $x$ be an intervention and $\hat{h}=\min(h^{(-)}(x),\ENDEPOCH)$.
If $\SMAXT{\hat{h}}(x)=C$ (the stock level goes to $C$ in the interval and before the first epoch with a stockout loss, if any), then, for all $\delta>0$ the following properties hold for intervention $x'= x+\delta$
\begin{align*}
    & L(x')=L(x)+\delta\\
    & \ENDSTOCK{x'}=\ENDSTOCK{x}
\end{align*}
where the new loss is of type surplus and is realized at one or more epochs earlier than $\hat{h}$.
\end{proposition}

\subsection{Optimal intervention characterization} \label{sec:OptimalityConditions}
Here we characterize the properties of optimal interventions for the uncapacitated problem and describe the set of optimal interventions. The computation of optimal interventions for the capacitated problem and the corresponding optimal value are then described.

\begin{proposition}\label{PN5:OptimalityConditions}
An intervention $x$ 
is optimal for the uncapacitated problem 
if and only if at least one of the following conditions holds:
\begin{align}
    & h^{(-)}(x)=h^{(+)}(x) = +\infty\label{PN5_09:OptimalityConditions_A}
    \\
    & h^{(+)}(x)<h^{(-)}(x) \quad \wedge \quad
        \SMINT{h^{(+)}(x)}(x) = 0
        \label{PN5_10:OptimalityConditions_BC}&
    \\
    & h^{(-)}(x)<h^{(+)}(x) \quad \wedge \quad
        \SMAXT{h^{(-)}(x)}(x) = C
        \label{PN5_11:OptimalityConditions_DE}&
\end{align}
\end{proposition}

The next proposition characterizes the set \XISSTAR\ of optimal interventions for the uncapacitated problem and its extreme points $\XIDSTAR=\min(\XISSTAR)$ and $\XIUSTAR=\max(\XISSTAR)$.

\begin{proposition}\label{PN6:MultipleOptimalInterventions}
Let $x^*$ be an optimal intervention for the uncapacitated problem. The whole set of optimal interventions is characterized as follows.
\begin{enumerate}
    \item If $L(x^*)>0$, 
    we have $\XISSTAR=\{x^*\}$, and for all $x'$
    \begin{align}\begin{array}{rl}\label{PN6_12:MOI.1}
        L(x') &=L(x^*)+|x'-x^*| \text{ and}\\
        \ENDSTOCK{x'}& = \ENDSTOCK{0}\ =\ENDSTOCK{x^*}.
    \end{array}\end{align}

    \item If $L(x^*)=0$, 
    the set of optimal interventions is a range \XIRSTAR\ defined as $$\XISSTAR=[x^*-\SMINT{\ENDEPOCH}(x^*),x^*+C-\SMAXT{\ENDEPOCH}(x^*)].$$

    Moreover, we have 
    \begin{align}\begin{array}{rll}\label{PN6_13:MOI.2}
    x'\in \XIRSTAR  \Rightarrow                         & L(x')=0,             & \ENDSTOCK{x'}=\ENDSTOCK{x^*}+(x'-x^*),\\
    x'>\XIUSTAR \Rightarrow & L(x')=x'-\XIUSTAR>0, & \ENDSTOCK{x'}=\ENDSTOCK{\XIUSTAR},\\
    x'<\XIDSTAR \Rightarrow & L(x')=\XIDSTAR-x'>0, & \ENDSTOCK{x'}=\ENDSTOCK{\XIDSTAR}.
    \end{array}
    \end{align}
\end{enumerate}
\end{proposition}

As a consequence of \PROREF{PN6:MultipleOptimalInterventions} we can conclude that the set of optimal interventions for the uncapacitated problem is an interval 
that
reduces to a single point when the optimal loss is positive.
In all cases any intervention $x\notin\XIRSTAR$ generates an additional loss equal to the distance of $x$ from the boundary of \XIRSTAR.

From now on, we denote by $\LISTAR$ the optimal loss and with $\XISTAR$ the optimal intervention with minimum modulus in \XIRSTAR, that is $\XISTAR=\XIDSTAR$ if $\XIDSTAR>0$, $\XISTAR=\XIUSTAR$ if $\XIUSTAR<0$, $\XISTAR=0$, otherwise.

The following proposition states the relationship between optimal interventions
for the uncapacitated problem and its capacitated version, together with the corresponding optimal values and stock level induced at the final epoch \ENDEPOCH.
Accordingly, we will denote by \XQSSTAR\ the set of optimal interventions for the capacitated problem with extreme points $\XQDSTAR=\min(\XQSSTAR)$ and $\XQUSTAR=\max(\XQSSTAR)$. We also denote by $\XQSTAR$ the optimal intervention with minimum modulus and as \LQSTAR\ the corresponding optimal value.

\begin{proposition}\label{PN7:XQfromXI}
The set \XQSSTAR\ of optimal interventions for the capacitated  version is a  range \XQRSTAR\ determined from the  range \XIRSTAR\ of optimal interventions for the uncapacitated \ONEINTPROBLEM\ problem as follows,
\begin{align}\label{PN7_14:PALGO.2.eqXR}
        \XQRSTAR = \left\{\begin{array}{ll}
        \{q\} & \text{if } q<\XIDSTAR=\XISTAR \\
        \{q-Q\} & \text{if } \XISTAR=\XIUSTAR<q-Q\\
        \XIRSTAR\cap [q-Q,q] & \text{otherwise,}
        \end{array}\right.
\end{align}
the minimum modulus optimal intervention \XQSTAR\ is determined from \XISTAR\ as
\begin{align}\label{PN7_15:PALGO.2.eq}
        \XQSTAR = \left\{\begin{array}{ll}
        q & \text{if } q<\XIDSTAR=\XISTAR \\
        q-Q & \text{if } \XISTAR=\XIUSTAR<q-Q\\
        \XISTAR & \text{if } \text{otherwise,}\\
        \end{array}\right.
\end{align}
the optimal loss is
\begin{align}\label{PN7_16:PALGO.4.eq3}
\LQSTAR =\LISTAR+|\XISTAR-\XQSTAR|,
\end{align}
and the stock level induced at epoch \ENDEPOCH\ is
\begin{align}\label{PN7_17:PALGO.4.eq4}
        \ENDSTOCK{\XQSTAR} = \ENDSTOCK{\XISTAR}=\ENDSTOCK{0}.
\end{align}
\end{proposition}

\subsection{Interval optimal intervention}\label{sec:intervaloptimization}

In this section we present an algorithm for the \ONEINTPROBLEM\ problem discussed above.
For the sake of generality and in view of how the algorithm will be used to solve the general problem,
we assume that the vectors $\SHBF$, $\lpBF$, $\lmBF$ and $\SBF$, representing, respectively, virtual stock, surplus loss, stockout loss and stock level induced by the null vector of interventions $\XBF=\OBF$ on all epochs $T$, have been already computed.

The idea is to maintain a $4$-tuple of variables
$\langle \XISTAR, \LISTAR, \DOWN, \UP\rangle$ where
\XISTAR\ represents the minimum modulus (optimal) intervention at  epoch \STARTEPOCH\ to minimize the total loss on interval $[\STARTEPOCH,h]$ for the uncapacitated problem;
variables \LISTAR, \DOWN\ and \UP, represent the optimal loss, the minimum and the maximum stock level induced by \XISTAR\ on the considered interval.
The computation is initialized by setting the $4$-tuple of variables for interval
$[\STARTEPOCH,\STARTEPOCH]$.
Then, the time interval is scanned to iteratively enlarge the considered interval up to
$[\STARTEPOCH,\ENDEPOCH]$.
Improvement opportunity for current intervention with respect to the new interval is evaluated, and by the end of each iteration the $4$-tuple is consistently updated to an interval with one more epoch.
When the scan is over, the optimal intervention \XQSTAR\ for the capacitated problem is computed together with the corresponding loss \LQSTAR.
A formal definition of the algorithm is provided in Pseudo-code \CiteAlgoVehicleIntervention\ (\ALGREF{A:\AlgoVehicleIntervention}) provided below.

\bigskip
\begin{algorithm}[H]\label{A:\AlgoVehicleIntervention} 
	\caption{\CiteAlgoVehicleIntervention}
	{
		\DontPrintSemicolon
		\SetKwInOut{Input}{input}\SetKwInOut{Output}{output}\SetKwInOut{Global}{global}
		\Input{
		$\STARTEPOCH:$ intervention epoch,
		$\ENDEPOCH:$ epoch up to which optimize intervention,
		$q:$ vehicle load,
		$Q:$ vehicle capacity
		}
		\Output{
		$\XQSTAR:$ minimum modulus optimal intervention for capacitated problem;
		$\LQSTAR:$ minimum loss for capacitated problem;
		$\XISTAR:$ minimum modulus optimal intervention for uncapacitated problem;
		$\LISTAR:$ minimum loss for uncapacitated problem;
		}
		\Global{$C$, $\SBF$, $\lmBF$, $\lpBF$}

        \smallskip
        \tcp*[l]{initialization}
        $\XISTAR = - l^{(\STARTEPOCH)+} + l^{(\STARTEPOCH)-}$\;\label{step:firstupdateX}
        $\LISTAR = 0$\;\label{step:firstupdateL}

        \smallskip
        $\DOWN, \UP = s^{(\STARTEPOCH)}$\;\label{step:firstupdateBOUNDS}

        \For(\tcp*[h]{scan time horizon}\label{1-int:mainloop}){$(h=\STARTEPOCH+1;\ h\leq \ENDEPOCH;\ h\texttt{++})$}{
            \tcp*[l]{intervention improvement}
            $\delta^{+} = \min(l^{(h)+},\DOWN)$\;\label{step:recoversurplus}
            $\delta^{-} = \min(l^{(h)-}, C-\UP)$\;\label{step:recoverstockout}

            \smallskip
            $\XISTAR = \XISTAR - \delta^{+} + \delta^{-}$\;\label{step:updateX}
            $\LISTAR = \LISTAR +(l^{(h)+} - \delta^{+}) + (l^{(h)-} - \delta^{-})$\;\label{step:updateL}

            \smallskip
            $\DOWN = \min(\DOWN - \delta^{+} ,s^{(h)})$\;\label{step:updateD}
            $\UP = \max(\UP + \delta^{-},s^{(h)})$\;\label{step:updateU}

        }
        \lIf{$\XISTAR<q-Q$}{
            $\{\XQSTAR=q-Q;\quad \LQSTAR=\LISTAR+\XQSTAR-\XISTAR\}$
        }\label{step:fix1}
        \lElseIf{$\XISTAR>q$}{
            $\{\XQSTAR=q;\quad \LQSTAR=\LISTAR+\XISTAR-\XQSTAR\}$
        }\label{step:fix2}
        \lElse{
            $\{\XQSTAR=\XISTAR;\quad \LQSTAR=\LISTAR\}$
        }\label{step:fix3}

    }
\end{algorithm}

\begin{theorem}\label{TH2:alg:Vehicle-Intervention}
Algorithm \CiteAlgoVehicleIntervention{} computes \XISTAR, \XQSTAR, \LISTAR\ and \LQSTAR\ in $O(\ENDEPOCH-\STARTEPOCH)$ time.
\end{theorem}

So far we have considered the initial stock level $\SZERO$ as given. However, the optimal interventions and losses for both the uncapacitated and capacitated problems can be somehow affected by the initial stock level.
In order to understand some properties of the general problem on the whole time horizon $T$, it can be useful to discuss a little bit about this point.
Let us denote by $\LI{s}{x}$ and $s^{(\ENDEPOCH)}(x\ |\ s)$ the total loss and the final stock level induced at epoch \ENDEPOCH\ by intervention $x$ given a stock level $s$ at epoch $\STARTEPOCH-1$.
Moreover, let us denote by $\XISTAR(s)$ and $\XQSTAR(s)$ the minimum modulus optimal interventions for the uncapacitated and capacitated problems, respectively, as functions of the stock level $s$ at epoch $\STARTEPOCH-1$,
and by $\LISTAR(s)$ and $\LQSTAR(s)$ the corresponding optimal losses.
Finally, we denote by $\XISSTAR(s)=[\XIDSTAR(s),\XIUSTAR(s)]$ and $\XQSSTAR(s)=[\XQDSTAR(s),\XQUSTAR(s)]$ the range of optimal interventions for the uncapacitated and capacitated problems, respectively.

The next proposition shows the relationship between these quantities when, under some circumstances, the initial stock level $s$ is translated by a quantity $\delta$.
The special cases considered by the proposition will be used in the following to study the properties of the general problem.

\begin{proposition}\label{PN8:PALGO.S0}
Let $s,s+\delta\in[0,C]$ be two alternative stock levels at epoch $\STARTEPOCH-1$, then
\begin{align}
   & \left\{\begin{array}{ll}
    \LISTAR(s+\delta)=\LISTAR(s)\DEFINE\LISTAR,\\
    \XISSTAR(s+\delta)=\XISSTAR(s)-\delta
    \end{array}\right.
    \label{PN8_18:PALGO.SO.1}\\
    & \LISTAR>0 & \Rightarrow &
    \quad s^{(\ENDEPOCH)}(x\ |\ s)=s^{(\ENDEPOCH)}(0\ |\ 0),
    \label{PN8_19:PALGO.SO.4}\\
    & \XISTAR(s)\geq q,\ \delta\leq\XISTAR(s)-q   & \Rightarrow &
    \left\{\begin{array}{ll}
    \XISTAR(s+\delta)=\XISTAR(s)-\delta\\
    \XQSTAR(s+\delta)=\XQSTAR(s)=q\\
    \LQSTAR(s+\delta)=\LQSTAR(s)-\delta\\
    s^{(\ENDEPOCH)}(\XQSTAR(s+\delta)\ |\ s+\delta)=s^{(\ENDEPOCH)}(\XQSTAR(s)\ |\ s),
    \end{array}\right.
    \label{PN8_20:PALGO.SO.2}\\
    & \XISTAR(s)\leq q-Q,\ \delta\geq\XISTAR(s)-(q-Q) & \Rightarrow &
    \left\{\begin{array}{ll}
    \XISTAR(s+\delta)=\XISTAR(s)-\delta\\
    \XQSTAR(s+\delta)=\XQSTAR(s)=q-Q\\
    \LQSTAR(s+\delta)=\LQSTAR(s)+\delta\\
    s^{(\ENDEPOCH)}(\XQSTAR(s+\delta)\ |\ s+\delta)=s^{(\ENDEPOCH)}(\XQSTAR(s)\ |\ s).
    \end{array}\right.
    \label{PN8_21:PALGO.SO.3}
\end{align}
\end{proposition}

Equations \EQUREF{PN8_18:PALGO.SO.1} show that the total loss induced in the interval by any intervention contains a ``systemic'' component which cannot be recovered in any way, independently of the initial stock level and available capacity and load.
On the other hand, the range of optimal interventions for the uncapacitated problem translate accordingly to the initial stock level; namely, an upshift of the initial stock level produces a corresponding downshift of the optimal interventions, and vice versa.
Also equation \EQUREF{PN8_19:PALGO.SO.4} shows a sort of stability property, when the systemic loss is positive any intervention with any initial stock level variation does not impact the final stock level.
The remaining part of the loss can be fully, or only partly, recovered according to the load-capacity constraint (i.e., we could recover more if we had more capacity or load) and to the initial stock level (we could lose less if we had a different initial stock).
While we have no margin to work on capacity, we can work on earlier interventions in order to set up appropriate initial stock levels for later interventions and even anticipate the recover of losses that could also be recovered later.
However, note that equations \EQUREF{PN8_20:PALGO.SO.2} and \EQUREF{PN8_21:PALGO.SO.3} show that, under applicable conditions, by changing the initial stock level we may improve or worsen the result according to the sign of stock change $\delta$.
This is the main focus of the analysis when we consider the optimization of the vector of interventions on the whole time horizon.

\section{Time horizon optimization} \label{sec:globalOpt}

In this section we first sketch the idea of how we can exploit the local optimization on the \ONEINTPROBLEM\ problem discussed in the above section to obtain the global optimization on the whole time horizon $T$ with $\NUMV$ intervention epochs.
Then, we define the algorithm and prove its correctness and computational complexity.

Let us
consider two families of time intervals made of $\NUMV$ members each
\begin{align*}\begin{array}{lll}
 I_i=[e_i,e_i+1,\ldots,e_{i+1}-1]             &\text{ for } i=1,\ldots,\NUMV-1; & \text{ and } I_\NUMV=[e_\NUMV,\ldots,m];\\
 J_i=[e_1,e_2,\ldots,e_{i+1}-1]             &\text{ for } i=1,\ldots,\NUMV-1; & \text{ and } J_\NUMV=[e_1,\ldots,m];
\end{array}\end{align*}
where $J_i=\cup_{l=1}^{i}I_l$ for $i=1,\ldots,\NUMV$ (recall we assumed w.l.o.g. $e_1=1$).

In the following, we denote by
$$\LVECTBKW{i} \DEFINE\sum_{h=e_1}^{e_{i+1}-1}\left(l^{(h)+}(\XBF) + l^{(h)-}(\XBF)\right)$$ 
the loss in interval $J_i$ depending on the initial stock level $\SZERO$ and the interventions $\XVECTrow{1}{i}$ implemented at epochs $e_1,\ldots,e_i$, and as
$$\LIOTT{i}\DEFINE\sum_{h=e_i}^{e_{i+1}-1}\left(l^{(h)+}(\XBF) + l^{(h)-}(\XBF)\right)$$
the loss in interval $I_i$ depending on the intervention $x_i$ implemented at epoch $e_i$ and on the stock level $s^{(e_i-1)}(\XVECTrow{1}{i-1})$ determined at epoch $e_i-1$ by all previous interventions and which we denote by \SX[x]{i}\ for short,
where for sake of notation compactness, when $i=1$, $\SX[x]{1}$  represents $\SZERO$.

By construction, we have the following equalities for all $i\leq\NUMV$
\begin{align}
\LVECTBKW{i} = \sum_{h=1}^{i}\LIOTT[x]{h}\label{eq22:EQ:sumLIOTT}
\end{align}
and
{\small
\begin{align}
& \begin{array}{rlr}
    \LVECTBKWONE & = \LIOTTGIVENS{1}{\SZERO} ,\\ 
    \LVECTBKW{i} & = \LVECTBKW{i-1}+\LIOTT[x]{i}, &\text{ for } i>1. 
\end{array}\label{eq23:EQ:LBCKWLIOTT}
\end{align}
}

\subsection{Backward formulation}

We see the computation of an optimal vector 
of interventions as a decision-making process in which we first determine the value of the first intervention at epoch $e_1$, then that of intervention at epoch $e_2$ and so on until we determine the last intervention at epoch $e_\NUMV$.
The choices made up to epoch $e_i$ condition the subsequent choices by influencing the stock level at epoch $e_{i+1}-1$ from which subsequent decisions take place.
Thus, once the first $i-1$ values of a vector $\XBF$ have been chosen, it will be a question of choosing the next one in order to minimize the losses in accordance with the conditions generated by the previous choices.
We will therefore have to choose at best the first $i$ values (i.e., $\XVECTrow[x]{1}{i}$) to obtain, after the optimization of the following ones (i.e., $\XVECTrow[x]{i+1}{\NUMV}$), a minimum overall loss.

Using equation \EQUREF{eq23:EQ:LBCKWLIOTT}, problem \EQUREF{eq07:base:L} can thus be formulated in a backward recursive form:
\begin{align}\label{eq24:base:Lrecursive:bkw}
    & \PROBLEMA{\NUMV}\DEFINE\min_{\XVECTrow[x]{1}{\NUMV-1}}\left(\LVECTBKW{\NUMV-1}+\min_{x_\NUMV\in[q_\NUMV-Q_\NUMV,q_\NUMV]}\LIOTT{\NUMV}\right)
\end{align}
where vector $\XVECTrow{1}{\NUMV-1}$ is required to satisfy feasibility constraints \EQUREF{eq01:update_q}.
Analogously, we define problems \PROBLEMA{i} for all $1<i\leq\NUMV$. For $i=1$, \PROBLEMA{i} reduces to the \ONEINTPROBLEM\ problem we have already solved in the previous section. We denote by \LSTAR\ the optimal value of \PROBLEMA{i} for any $i\leq\NUMV$.

According to formulation \EQUREF{eq24:base:Lrecursive:bkw}, for any given vector $\XVECTrow[x]{1}{\NUMV-1}$ the best we can do is to solve the \ONEINTPROBLEM\ problem on interval $I_\NUMV$ given the initial stock $\SX{\NUMV}$ produced at epoch $e_\NUMV-1$ by vector $\XVECTrow[x]{1}{\NUMV-1}$.
Thus, we have to find the best vector $\XVECTrow[x]{1}{\NUMV-1}$ such that after optimizing $x_\NUMV$ we get the overall minimum loss.

\bigskip
We may observe that if $\SX{\NUMV}$ were independent of $\XVECTrow[x]{1}{\NUMV-1}$, that is $\SX{\NUMV}=S_\NUMV$ for a fixed and known $S_\NUMV$, then we could solve separately the two (smaller) problems:
first the \ONEINTPROBLEM\ problem
$$
\min_{x_\NUMV\in[Q_\NUMV-q_\NUMV,q_\NUMV]}\LIOTTGIVENS{\NUMV}{S_\NUMV}
$$
and then
$$
\PROBLEMA{\NUMV-1}=\PIOTTBKW{\NUMV-1}
$$
which could be recursively broken in $\NUMV-1$ independent \ONEINTPROBLEM\ problems.
Although this is true in some cases, it does not hold in general
because the interventions on $J_{\NUMV-1}$ may affect, for better or worse, the initial stock $S_w$ for $I_{\NUMV}$.
Nevertheless, our approach to problem \EQUREF{eq24:base:Lrecursive:bkw} is to prove that we can solve two separate smaller problems of these types leading to the solution of a sequence of \ONEINTPROBLEM\ problems.
The following discussion illustrates how we can separate the two subproblems.

\begin{proposition}\label{PN9:OttimoOgniIntervallo}
Let $i\leq\NUMV$ and $\XBF^*=\XVECTrow[x^*]{1}{i}$ be an optimal vector of interventions for problem \PROBLEMA{i}.
Then, for all $h\leq i$, intervention $x^*_h$ is optimal for the \ONEINTPROBLEM\ capacitated problem on interval $I_h$ given the stock level $\SX[x^*]{h}$ at epoch $e_{h}-1$. That is,
$$\LIOTT[x^*]{h}=\underset{x\in[q_h-Q_h,q_h]}{\min}\ \LIXS[x]{h}{\SX[x^*]{h}}.$$
\end{proposition}

\begin{proposition}\label{PN10:ottimoIntervalli} 
    Let $\XBF^*$ and $\YBF^*$ be two optimal vectors of interventions for \PROBLEMA{i}.
    Then, for any interval $I_h$ ($h\leq i$) 
    we have
        \begin{align}
        & \LIOTT[x^*]{h}=\LIOTT[y^*]{h}.\label{PN10_25:EQ:14.UGL}
        \end{align}
\end{proposition}

\begin{proposition}\label{PN11:OttimoTotaleParziale}
Let $\XBF^*=\XVECTrow[x^*]{1}{i}$ be an optimal vector of interventions for problem \PROBLEMA{i} for some $1<i\leq\NUMV$. Then for all $h<i$, vector $\XVECTrow[x^*]{1}{h}$ is optimal for
\PROBLEMA{h}.
\end{proposition}

Observe that while \PROREF{PN11:OttimoTotaleParziale} cannot guarantee that an optimal vector of interventions for problem \PROBLEMA{\NUMV} can be built upon any optimal intervention for problem \PROBLEMA{i} with $i<\NUMV$, it however guarantees that in order to find the first $i$ optimal interventions for problem \PROBLEMA{\NUMV} we only have to look into the set of optimal vectors of interventions for problem \PROBLEMA{i}.

\begin{proposition}\label{PN12:neutralStockOptimum}
Problem \PROBLEMA{i}\
admits an optimal vector of interventions $\XBF^*$ such that $s^{(f)}(\XBF^*)=s^{(f)}(\OBF)$ where $f=e_{i+1}-1$ if $i<\NUMV$, $f=m$ otherwise.
\end{proposition}

\subsection{The augmented problem}

Let us consider, for $i<\NUMV$ and a scalar parameter $\delta$, the augmented time intervals
$\AUGMINTI=[e_i,e_{i+1}]$, and $\AUGMINTJ=[e_1,e_{i+1}]$
(i.e., $\AUGMINTI=I_i\cup\{e_{i+1}\}$ and $\AUGMINTJ=J_i\cup\{e_{i+1}\}$) where net flow $d^{(e_{i+1})}$ is fictitiously set to zero if $\delta=0$ or so that the null vector $\OBF$ produces, at epoch $e_{i+1}$ a stockout loss $|\delta|$ if $\delta>0$ or a surplus loss $|\delta|$ if $\delta<0$.
That is, if $\delta=0$, $d^{(e_{i+1})}=0$, if $\delta > 0$, $d^{(e_{i+1})}=-s^{(e_{i+1}-1)}(\OBF)-|\delta|$, and, if $\delta<0$, $d^{(e_{i+1})}=C-s^{(e_{i+1}-1)}(\OBF)+|\delta|$.

We denote by $\LVECTBKWAUGMN{i}$ the corresponding loss induced by interventions $\XVECTrow[x]{1}{i}$ and call
\begin{align*}
    \PROBLEMAAUG{i}\DEFINE\PIOTTBKWAUGMN{i}
\end{align*}
the \textit{augmented problem} on interval $\AUGMINTJ$. We generally denote by \LSTARAUGMN[i] the optimal value of \PROBLEMAAUG{i}.

Note that the augmented problem is still a problem in the class defined by \EQUREF{eq07:base:L} and all derived properties still hold for it.
However, \PROBLEMA{i} and \PROBLEMAAUG{i} are defined on different set of data; namely problem \PROBLEMAAUG{i} has the same set of decision epoch as \PROBLEMA{i}, but a longer time horizon with one more epoch: $e_{i+1}$ and net flow $d^{(e_{i+1})}$.
 Finally, we have by construction
             \begin{align*}
             \LVECTBKWOAUGMN{i} & = \LVECTBKWO{i} + |\delta|\\
             \LIOTTAUGMNX[0]{i}{\SX[0]{i}} & = \LIXS[0]{i}{\SX[0]{i}} + |\delta|\\
             \end{align*}

We can further state some more properties of the augmented problem \PROBLEMAAUG{i} that establish its relationship with the base problem \PROBLEMA{i}.

\begin{proposition}\label{PN13:augmentation}
For all $i<\NUMV$, an optimal vector of interventions \XSTARAUGMN\ for the augmented problem \PROBLEMAAUG{i} is optimal also for the problem \PROBLEMA{i}. In a formula:
$$\LVECTBKW[\bar{x}]{i}=\PIOTTBKW{i}=\LSTAR.$$
\end{proposition}

\begin{proposition}\label{PN14:MERGE} 
    Let us consider problem \PROBLEMA{i} for some $i>1$ and
    let
    $S_i=\INISTOCKEPOCHBYVECT[\OBF]{i}$ be the stock induced by the null vector of interventions at epoch $e_{i}-1$.
    Let
    $(\XISTAR(S_i),\LISTAR_{I_i})$ and $(\XQSTAR(S_i),\LQSTAR_{I_i}(S_i))$
    be the pairs of optimal intervention and optimal value for
    the uncapacitated and capacitated \ONEINTPROBLEM\ problems
    on interval $I_i$ with initial stock level $S_i$.
    Then, an optimal vector of interventions 
    $\XSTARAUGMN=\XVECTrow[\bar{x}]{1}{i-1}$ for the augmented problem \PROBLEMAAUG{i-1}, with $\delta=\XISTAR(S_i)-\XQSTAR(S_i)$, exists such that 
    $\XBF^*=[\bar{x}_1,\ldots,\bar{x}_{i-1},\XQSTAR(S_i)]$ is an optimal vector of interventions for problem \PROBLEMA{i} and
    $\LSTAR[i]=\LSTARAUGMN[i-1]+\LISTAR$.
\end{proposition}

\subsection{The algorithm}

\subsubsection{Sketch of the algorithm}

The algorithm exploits the backward formulation of the problem. The basic line of the algorithm is described below.
We first compute the vectors $\SHBF$, $\lpBF$, $\lmBF$, and $\SBF$ representing, respectively, virtual stock, surplus loss, stockout loss and stock level induced by the null vector of interventions $\XBF=\OBF$ on all epochs $T$.
Then a \ONEINTPROBLEM\ problem is solved on the last interval $I_\NUMV$ with initial stock level $S_\NUMV=s^{(e_\NUMV-1)}(\OBF)$. The difference $\XISTAR_{I_\NUMV}-\XQSTAR_{I_\NUMV}$ is used to build the augmented problem on interval $J_{\NUMV-1}$. In some sense we may say that $|\XISTAR_{I_\NUMV}-\XQSTAR_{I_\NUMV}|$ is the loss due to limited capacity of vehicle $v_\NUMV$ and which is delegated to previous interventions.
The augmented problem is recursively solved by solving a \ONEINTPROBLEM\ problem on the last interval and delegating to previous decisions the loss that could not be recovered by last intervention due to limited capacity.
Recursion stops when the augmented problem is defined on the first interval; in that case only the \ONEINTPROBLEM\ problem is solved.
In the meanwhile we take into account the loss realized on each interval and compose the vector of optimal interventions.

\subsubsection{Algorithm definition}

A formal definition of the algorithm is provided in Pseudo-codes \CiteAlgoGlobalBackward{} (\ALGREF{A:\AlgoGlobalBackward}) and \CiteAlgoGlobalBackwardRec{} (\ALGREF{A:\AlgoGlobalBackwardRec}) provided below.

\begin{algorithm}[H]\label{A:\AlgoGlobalBackward} 
	\caption{\CiteAlgoGlobalBackward}
	{
		\DontPrintSemicolon
		\SetKwInOut{Input}{input}\SetKwInOut{Output}{output}\SetKwInOut{Global}{global}
		\Input{$s:$ initial stock level, $C:$ \STATION\ capacity, $\DBF:$ net flow, $H:$ decision epochs, $\qBF:$ vehicle loads, $\QBF:$ vehicle capacities}
		\Output{$\XBF^*$ optimal vector of interventions; $L^*$ optimal value}
		\Global{all variables}

        \medskip
        \tcp*[l]{initialization}
        $m=\texttt{dim}(\DBF)$ \tcp*[l]{problem size}\label{TH3_1}
        $\NUMV=|H|$ \tcp*[l]{number of interventions}
        $\SZERO = s$\;
        \For{$(h=1:\ h\leq m;\ h\texttt{++})$}{
            $\hat{s}^{(h)} = s^{(h-1)}+d^{(h)}$ \tcp*[l]{$\SHBF=\SHBF(\OBF)$}
            $l^{(h)+} = \max(0,\hat{s}^{(h)}-C)$ \tcp*[l]{$\lpBF=\lpBF(\OBF)$}
            $l^{(h)-} = \max(0,-\hat{s}^{(h)})$ \tcp*[l]{$\lmBF=\lmBF(\OBF)$}
            $s^{(h)} = \hat{s}^{(h)} - l^{(h)+} + l^{(h)-}$ \tcp*[l]{$\SBF=\SBF(\OBF)$}\label{TH3_8}
        }

        \medskip
        $\STARTEPOCH=e_w$\;\label{TH3_9}
        $\ENDEPOCH = m$\;\label{TH3_10}
        $(\XQSTAR,\LQSTAR,\XISTAR,\LISTAR)=\CiteAlgoVehicleIntervention(\STARTEPOCH,\ \ENDEPOCH,\ q_w,\ Q_w)$\;\label{TH3_11}
        $(\XBF,L)=\CiteAlgoGlobalBackwardRec(w-1,\XISTAR - \XQSTAR)$\;\label{TH3_12}
        $\XBF^*=[\XBF, \XQSTAR]$\;\label{TH3_13}
        $L^*=L + \LISTAR$\;\label{TH3_14}

        }
\end{algorithm}

\begin{algorithm}[H]\label{A:\AlgoGlobalBackwardRec}
	\caption{\CiteAlgoGlobalBackwardRec}
	{
		\DontPrintSemicolon
		\SetKwInOut{Input}{input}\SetKwInOut{Output}{output}\SetKwInOut{Global}{global}
		\Input{$i:$ stage of decision process, $\delta:$ problem augmentation}
		\Output{$\XBF^*$ optimal vector of interventions up to stage $i$; $L^*$ optimal value up to stage $i$}
		\Global{$C$, $\DBF$, $\SBF$, $H$, $\qBF$, $\QBF$}

        \medskip
        \If{$i>0$}{\label{TH3R_1}
            $\STARTEPOCH=e_i$\;\label{TH3R_2}
            $\ENDEPOCH = e_{i+1}$\;\label{TH3R_3}
            \lIf {$\delta> 0$}{$d^{(\ENDEPOCH)}=-s^{(\ENDEPOCH-1)}-|\delta|$}\label{TH3R_4}
            \lElseIf{$\delta< 0$}{$d^{(\ENDEPOCH)}=C-s^{(\ENDEPOCH-1)}+|\delta|$}\label{TH3R_5}
           \lElse{$d^{(\ENDEPOCH)}=0$}\label{TH3R_5a}
            $(\XQSTAR,\LQSTAR,\XISTAR,\LISTAR)=\CiteAlgoVehicleIntervention(\STARTEPOCH,\ \ENDEPOCH,\ q_i,\ Q_i)$\;\label{TH3R_6}
            $(\XBF,L)=\CiteAlgoGlobalBackwardRec(i-1,\XISTAR - \XQSTAR)$\;\label{TH3R_7}
            $\XBF^*=[\XBF, \XQSTAR]$\;\label{TH3R_8}
            $L^*=L+\LISTAR$\;\label{TH3R_9}
        }
        \Else{\label{TH3R_10}
            $\XBF^*=\emptyset$\;\label{TH3R_11}
            $L^*=|\delta|$\;\label{TH3R_12}
        }
    }
\end{algorithm}

\subsubsection{Algorithm properties}

\begin{theorem}\label{PN:algoBKWcorrectness}
Algorithm \CiteAlgoGlobalBackward\ (\ALGREF{A:\AlgoGlobalBackward}) computes an optimal vector of interventions $\XBF^*$ for problem \PROBLEMA{\NUMV} in $O(m)$ time. 
Moreover, the vector returned by the algorithm  induces at the final epoch $m$ the same stock induced by the null vector $\OBF$ (i.e., $s^{(m)}(\XBF^*)=s^{(m)}(\OBF)$).
\end{theorem}

\section{Conclusions}

In this paper, the one-station bike repositioning problem is studied. A set of capacitated vehicles with a given bike load are planned to visit a station to load or unload bikes at given times of a finite time horizon. The problem is to find the optimal number of bikes to load or unload at each visit with the goal of minimizing the number of user requests (of bike rental or return) lost over the time horizon.
An optimal algorithm with linear complexity in the cardinality of the time horizon is presented. In addition to finding the optimal solution to the problem, the algorithm provides the number of requests that would be lost regardless of the capacity and initial load of the vehicles.

The studied problem considers multiple visits of vehicles to the same station and takes into account the interaction among vehicles, crucial aspects of the problem of repositioning bikes in a bike sharing system that have been in most of the cases disregarded in the literature. The algorithm and the results presented in this paper are a contribution to a deep understanding of the general repositioning problem and could be useful to take decisions, for example about the initial load of the vehicles. The algorithm, in particular, can be seen as a component of a general solution approach to a more general problem.

As a final remark, we would like to point out that the
strategy adopted by the proposed algorithm is of the kind ``recover as late as possible''.
However, different optimal repositioning strategies, like ``recover as soon as possible'', can be implemented with little changes to the algorithm, maintaining its linear complexity.

Future research efforts should be devoted to extending the analysis to the case where the times of the visits are not given and to the case where more than one station is considered.

\section*{Acknowledgements}
The authors wish to express their gratitude to the anonymous referees who provided very useful and detailed comments 
to improve the manuscript.

\appendix\appendixpage

\section{Proofs}\label{app:proof}

\subsection{\ONEINTPROBLEM\ problem}

\begin{myproof}[Proof of \PROREF{PN1:ImproveByDownShift}]

Let us iteratively build vectors $\SBF(x')$, $\lmBF(x')$ and $\lpBF(x')$.
It is easy to see that from equations \EQUREF{eq02:update_shat_int}-\EQUREF{eq06:update_s} we have the following cases.

\medskip
\noindent\textbf{Case $\STARTEPOCH=h^{(+)}(x)$.}

\begin{itemize}
\item For $k=\STARTEPOCH=h^{(+)}(x)$:
\end{itemize}
\begin{align*}
    \left\{
    \begin{array}{rl}
    \hat{s}^{(k)}(x')&=s^{(k-1)}+d^{(k)}+x'= s^{(k-1)}+d^{(k)}+(x-\delta)=\hat{s}^{(k)}(x)-\delta=C+l^{(k)+}(x)-\delta\geq C\\
    \BOXED{l^{(k)+}(x')}&=\max(0,\hat{s}^{(k)}(x')-C) 
                =\max(0,l^{(k)+}(x)-\delta) =\BOXED{l^{(k)+}(x)-\delta}\\
    \BOXED{l^{(k)-}(x')}&=\max(0,-\hat{s}^{(k)}(x')) = 0 = \BOXED{l^{(k)-}(x)}\\
    \BOXED{s^{(k)}(x')} &=\hat{s}^{(k)}(x')-l^{(k)+}(x')+l^{(k)-}(x') = s^{(k)}(x)+l^{(k)+}(x)-\delta -(l^{(k)+}(x)-\delta) \BOXED{= s^{(k)}(x)=C}.
    \end{array}
    \right.
\end{align*}

\medskip

\noindent\textbf{Case $\STARTEPOCH<h^{(+)}(x)$.}

\begin{itemize}
\item For $k=\STARTEPOCH$, we have:
\end{itemize}
\begin{align*}
    \left\{
    \begin{array}{rl}
    \hat{s}^{(k)}(x')&=s^{(k-1)}+d^{(k)}+x'=s^{(k-1)}+d^{(k)}+(x-\delta)=\hat{s}^{(k)}(x)-\delta\leq C\\%=s^{(k)}(x)-\delta\\%\in[0,C]\\
    \BOXED{l^{(k)+}(x')}&=\max(0,\hat{s}^{(k)}(x')-C) =0= \BOXED{l^{(k)+}(x)}\\
    \BOXED{l^{(k)-}(x')}&=\max(0,-\hat{s}^{(k)}(x')) = \max(0,-(\hat{s}^{(k)}(x)-\delta))=
    \BOXED{l^{(k)-}(x)}\\
    \BOXED{s^{(k)}(x')} &=\hat{s}^{(k)}(x')-l^{(k)+}(x')+l^{(k)-}(x')=\hat{s}^{(k)}(x)-l^{(k)+}(x')+l^{(k)-}(x')-\delta =\BOXED{s^{(k)}(x)-\delta}.
    \end{array}
    \right.
\end{align*}

\begin{itemize}
\item For $k=\STARTEPOCH+1,\ldots,h^{(+)}(x)-1$, we have:
\end{itemize}
\begin{align*}
    \left\{
    \begin{array}{ll}
    \hat{s}^{(k)}(x')&=s^{(k-1)}(x')+d^{(k)}= (s^{(k-1)}(x)-\delta)+d^{(k)}=(s^{(k-1)}(x)+d^{(k)})-\delta=
                \hat{s}^{(k)}(x)-\delta\leq C\\
    \BOXED{l^{(k)+}(x')} &=\max(0,\hat{s}^{(k)}(x')-C) = 0= \BOXED{l^{(k)+}(x)}\\
    \BOXED{l^{(k)-}(x')} &=\max(0,-\hat{s}^{(k)}(x')) = \max(0,-(\hat{s}^{(k)}(x)-\delta))=\BOXED{l^{(k)-}(x)}\\
    \BOXED{s^{(k)}(x')} &=\hat{s}^{(k)}(x')-l^{(k)+}(x')+l^{(k)-}(x') =\hat{s}^{(k)}(x)-l^{(k)+}(x')+l^{(k)-}(x')-\delta= \BOXED{s^{(k)}(x)-\delta}.
    \end{array}
    \right.
\end{align*}

\begin{itemize}
\item For $k=h^{(+)}(x)$, we have:
\end{itemize}
\begin{align*}
    \left\{
    \begin{array}{rl}
    \hat{s}^{(k)}(x') &=s^{(k-1)}(x')+d^{(k)}= (s^{(k-1)}(x)-\delta)+d^{(k)}=\hat{s}^{(k)}(x)-\delta=C+l^{(k)+}(x)-\delta\geq C \\
    \BOXED{l^{(k)+}(x')} &=\max(0,\hat{s}^{(k)}(x')-C)=\max(0,\hat{s}^{(k)}(x)-C-\delta) =\max(0,l^{(k)+}(x)-\delta) = \BOXED{l^{(k)+}(x)-\delta}\\
    \BOXED{l^{(k)-}(x')} &=\max(0,-\hat{s}^{(k)}(x')) = \max(0,-(\hat{s}^{(k)}(x)-\delta))= 0=\BOXED{l^{(k)-}(x)}\\
    \BOXED{s^{(k)}(x')} &=\hat{s}^{(k)}(x')-l^{(k)+}(x')+l^{(k)-}(x') =s^{(k)}(x)+l^{(k)+}-\delta-(l^{(k)+}(x)-\delta)\BOXED{=s^{(k)}(x)= C}.
    \end{array}
    \right.
\end{align*}

\medskip
\noindent
Observe that interventions $x$ and $x'$ produce the same \STATION\ stock level $C$ at epoch $h^{(+)}(x)$. Thus, since then they produce the same virtual stock, losses and stock level up to \ENDEPOCH.

\medskip
As far as total loss is concerned, equation $L(x')=L(x)-\delta$ comes from the construction of $\lpBF(x')$ and clearly the recovered loss is equal to intervention increment.
Finally, we have  $h^{(+)}(x')>h^{(+)}(x)$ if and only if the whole surplus loss at epoch $h^{(+)}(x)$ is fully recovered, which happens if and only if $\delta$ equals the loss produced at that epoch by intervention $x$.
\end{myproof}

\bigskip
\begin{myproof}[Proof of \PROREF{PN2:ImproveByUpShift}]
Discussion mirroring that for \PROREF{PN1:ImproveByDownShift}.
\end{myproof}

\bigskip

\begin{myproof}[Proof of \PROREF{PN3:WorsenByDownShift}]
Let us consider the first epoch $\HZEROX$ such that $s^{(\HZEROX)}(x)=0$ and observe that the inequality $\HZEROX\leq \hat{h}$ holds by the hypothesis.

\noindent\textbf{Case $\HZEROX=\STARTEPOCH$}.

We have $\hat{s}^{(\STARTEPOCH)}(x)\leq0$ and every $\delta>0$ implies
$\hat{s}^{(\STARTEPOCH)}(x')=\hat{s}^{(\STARTEPOCH)}(x)-\delta<0$
with $l^{(\STARTEPOCH)-}(x')=l^{(\STARTEPOCH)-}(x)+\delta$.
Moreover, the equality $s^{(\STARTEPOCH)}(x')=s^{(\STARTEPOCH)}(x)=0$ holds, which implies all future values of $\SBF$, $\lpBF$ and $\lmBF$ will be the same for $x$ and $x'$; in particular, with $s^{(\ENDEPOCH)}(x')=s^{(\ENDEPOCH)}(x)$ and $L(x')=L(x)+\delta$.

Observe that a new stockout loss is introduced at epoch $\HZEROX$ while any other value remains unchanged.

\smallskip
\noindent\textbf{Case $\HZEROX>\STARTEPOCH$}.

For ease of notation, let us indicate as $\SMINBEFOREHZEROX=\SMINT{\HZEROX-1}(x)>0$ the minimum level reached by the \STATION\ stock before epoch $\HZEROX$.

\begin{itemize}
\item For $\delta\leq\SMINBEFOREHZEROX$ we can prove by construction that
$$\hat{s}^{(k)}(x')=\hat{s}^{(k)}(x)-\delta\in[0,C) \text{ for all } k\in[\STARTEPOCH,\HZEROX-1],$$
while
$$\hat{s}^{(k)}(x')=\hat{s}^{(k)}(x)-\delta \text{ for } k=\HZEROX$$
from which
$$l^{(\HZEROX)-}(x')=-\hat{s}^{(\HZEROX)}(x')=-(\hat{s}^{(\HZEROX)}(x)-\delta)=l^{(\HZEROX)-}(x)+\delta.$$
Moreover, the equality $s^{(\HZEROX)}(x')=s^{(\HZEROX)}(x)=0$ holds, which implies all future values of $\SBF$, $\lpBF$ and $\lmBF$ will be the same for $x$ and $x'$; in particular, with $s^{(\ENDEPOCH)}(x')=s^{(\ENDEPOCH)}(x)$ and $L(x')=L(x)+\delta$.

Observe that the stock is uniformly downshifted by $\delta$ up to $\HZEROX$ and a new stockout loss is introduced at epoch $\HZEROX$ while any other value remains unchanged.

\item For $\delta>\SMINBEFOREHZEROX$ we proceed as follows.
Let us define $\DBFI[1]=\SMINBEFOREHZEROX$ and $\XBFI[1]=x-\DBFI[1]$ the worsening intervention with $s^{(\ENDEPOCH)}(\XBFI[1])=s^{(\ENDEPOCH)}(x)$ and $L(\XBFI[1])=L(x)+\DBFI[1]$.
Observe that with respect to $x$, intervention \XBFI[1] uniformly downshifts the stock by \DBFI[1] up to $\HZEROX$ and introduces a new stockout loss \DBFI[1] at epoch $\HZEROX$ while any other value remains unchanged. Moreover, intervention $\XBFI[1]$ satisfies hypothesis of \PROREF{PN3:WorsenByDownShift}, and thus, we can iterate on $\XBFI[1]$ with $\delta'=\delta-\DBFI[1]$.

If $\delta'> \underline{s}'(\XBFI[1])$, we iterate in this case with $\DBFI[2]=\underline{s}'(\XBFI[1])$ and $\XBFI[2]=\XBFI[1]-\DBFI[2]$ producing a worse loss $L(\XBFI[2])=L(\XBFI[1])+\DBFI[2]=L(x)+\DBFI[1]+\DBFI[2]$ and so on.

We reach the last iteration when $\delta'\leq \underline{s}'(\XBFI[j])$, so that we fall in the first case and end up with
intervention $x'=\XBFI[j]-\delta'=x-\sum_{k=1}^{j}\DBFI[j])-\delta'=x-\delta$ producing a loss $L(x')=L(\XBFI[j])+\delta'=L(x)+\sum_{k=1}^{j}\DBFI[j]+\delta'=L(x)+\delta$.

Observe that at each iteration a new earlier stockout loss is introduced.
\end{itemize}
\end{myproof}

\bigskip
\begin{myproof}[Proof of \PROREF{PN4:WorsenByUpShift}]
Discussion mirroring that for \PROREF{PN3:WorsenByDownShift}.
\end{myproof}

\bigskip
\begin{myproof}[Proof of \PROREF{PN5:OptimalityConditions}]

    Sufficient conditions

    \begin{enumerate}
        \item[Eq. \EQUREF{PN5_09:OptimalityConditions_A}]
        If $h^{(+)}(x) = h^{(-)}(x) = +\infty$, then $L(x)=0$ holds by definition, which cannot be further improved.
        \item[Eq. \EQUREF{PN5_10:OptimalityConditions_BC}]
        Assuming $h^{(+)}(x)<h^{(-)}(x)$ (first loss is of type surplus),
        hypothesis of \PROREF{PN4:WorsenByUpShift} are satisfied (stock level reaches $C$ before any stockout loss). Thus any upshift brings to a larger loss.
        On the other hand, any downshift also leads to a worse loss. Indeed, $\SMINT{h^{(+)}(x)}(x)=0$, and, according to \PROREF{PN3:WorsenByDownShift}, any downshift would imply additional stockout loss at epochs earlier than  $h^{(+)}(x)$.
        \item[Eq. \EQUREF{PN5_11:OptimalityConditions_DE}]
        Discussion mirroring that for point \EQUREF{PN5_10:OptimalityConditions_BC}.
    \end{enumerate}

Necessary conditions. We proceed by contradiction. Let us assume that none of the listed conditions hold. Then it must be the case that $\min(h^{(+)}(x),h^{(-)}(x))<\infty$ (some kind of loss occurs at some epoch). Of course it cannot be the case that $h^{(+)}(x)=h^{(-)}(x)<\infty$, thus, exactly one of the following two holds

\begin{itemize}
    \item $h^{(+)}(x)<h^{(-)}(x)$. In this case we must have $\SMINT{h^{(+)}(x)}(x)>0$. Then hypothesis of \PROREF{PN1:ImproveByDownShift} hold with $\delta^*>0$ and losses produced by intervention $x$ can be strictly improved.
    \item $h^{(-)}(x)<h^{(+)}(x)$. In this case we have $\SMAXT{h^{(-)}(x)}(x)<C$. Then hypothesis of \PROREF{PN2:ImproveByUpShift} hold with $\delta^*>0$ and losses produced by intervention $x$ can be strictly improved.
\end{itemize}
\end{myproof}

\bigskip
\begin{myproof}[Proof of \PROREF{PN6:MultipleOptimalInterventions}]
Let $x^*$ be an optimal intervention and consider the cases $L(x^*)>0$ and $L(x^*)=0$ separately.
\begin{itemize}
    \item [$L(x^*)>0$.] At least one loss is detected in $[\STARTEPOCH,\ENDEPOCH]$ and either condition \EQUREF{PN5_10:OptimalityConditions_BC} or \EQUREF{PN5_11:OptimalityConditions_DE} is satisfied by $x^*$.
    Let us assume w.l.o.g. that condition \EQUREF{PN5_10:OptimalityConditions_BC} is matched by intervention $x^*$ (i.e., $h^{(+)}(x^*)<h^{(-)}(x^*)$ and $\SMINT{h^{(+)}(x^*)}(x^*)=0$) and consider an alternative intervention $x'\neq x^*$. We have to discuss two cases
    \begin{itemize}
        \item Case $x'=x^*-\delta$ with $\delta>0$.
        Condition \EQUREF{PN5_10:OptimalityConditions_BC} guarantees that hypothesis of \PROREF{PN3:WorsenByDownShift} are satisfied by $x^*$. Accordingly, any alternative intervention $x'<x^*$ would lead to a loss $L(x')=L(x^*)+|x'-x^*|$ with $s^{(\ENDEPOCH)}(x')=s^{(\ENDEPOCH)}(x^*)$.

        \item Case $x'=x^*+\delta$ with $\delta>0$.
        Since we have $h^{(+)}(x^*)<\infty$, we have a surplus loss in the interval with $\SMAXT{h}(x^*)=C$ for some $h\leq\ENDEPOCH$. Thus, hypothesis of \PROREF{PN4:WorsenByUpShift} are satisfied by $x^*$. Accordingly, any alternative intervention $x'>x^*$ would lead to a loss $L(x')=L(x^*)+|x'-x^*|$ with $s^{(\ENDEPOCH)}(x')=s^{(\ENDEPOCH)}(x^*)$.

    \end{itemize}

    \item[$L(x^*)=0$.] Condition \EQUREF{PN5_09:OptimalityConditions_A} is satisfied by $x^*$.
    First, observe that
    \begin{align*}
    \hat{s}^{(\STARTEPOCH)}(x')=s^{(\STARTEPOCH-1)}+d^{(\STARTEPOCH)}+x'=s^{(\STARTEPOCH-1)}+d^{(\STARTEPOCH)}+x^*+(x'-x^*)=\hat{s}^{(\STARTEPOCH)}(x^*)+(x'-x^*)
    \end{align*}
     with $\hat{s}^{(\STARTEPOCH)}(x^*)\in[0,C]$ by the hypothesis $L(x^*)=0$ and thus $s^{(\STARTEPOCH)}(x^*)=\hat{s}^{(\STARTEPOCH)}(x^*)$.

\medskip

    If $x'\in \XISSTAR$ we have $(x'-x^*)\in[-\SMINT{\ENDEPOCH}(x^*),C-\SMAXT{\ENDEPOCH}(x^*)]$ and
    thus $\hat{s}^{(\STARTEPOCH)}(x')=s^{(\STARTEPOCH)}(x^*)+(x'-x^*)\in[0,C]$, and also $s^{(\STARTEPOCH)}(x')=s^{(\STARTEPOCH)}(x^*)+(x'-x^*)\in[0,C]$.

    Moreover, for $h=\STARTEPOCH+1,\ldots,\ENDEPOCH$ we get
    \begin{align*}
    \hat{s}^{(h)}(x')=s^{(h-1)}(x')+d^{(h)}=s^{(h-1)}(x^*)+(x'-x^*)+d^{(h)}=\hat{s}^{(h)}(x^*)+(x'-x^*) 
    \end{align*}
    with $\hat{s}^{(h)}(x^*)\in[0,C]$ by the hypothesis $L(x^*)=0$ and thus $\hat{s}^{(h)}(x')\in[0,C]$ and also $s^{(h)}(x')=s^{(h)}(x^*)+(x'-x^*)$.

    This proves $L(x')=0$ and optimality of $x'$.

    Now we show that extreme points of $\XISSTAR$ do not depend on the particular optimal intervention $x^*$.
    Let us first observe that for any $x'\in \XISSTAR$, from $s^{(h)}(x')=s^{(h)}(x^*)+(x'-x^*)$ for all $h\in[\STARTEPOCH,\ENDEPOCH]$, we derive $\SMINT{\ENDEPOCH}(x')=\SMINT{\ENDEPOCH}(x^*)+(x'-x^*)$ and $\SMAXT{\ENDEPOCH}(x')=\SMAXT{\ENDEPOCH}(x^*)+(x'-x^*)$.
    Algebra tells us that $x^*-\SMINT{\ENDEPOCH}(x^*)=x'-\SMINT{\ENDEPOCH}(x')$ and $x^*+C-\SMAXT{\ENDEPOCH}(x^*)=x'+C-\SMAXT{\ENDEPOCH}(x')$ which is what we need.

\medskip
    If $x'>\XIUSTAR=\max(\XISSTAR)$ we observe that for the optimal intervention $\XIUSTAR$ we have $\SMAXT{\ENDEPOCH}(\XIUSTAR)=C$ and $\XIUSTAR$ satisfies hypothesis of \PROREF{PN4:WorsenByUpShift}. Thus if we write $x'=\XIUSTAR+(x'-\XIUSTAR)$ we get $L(x')=x'-\XIUSTAR$ and $\ENDSTOCK{x'}=\ENDSTOCK{\XIUSTAR}$.

\medskip
    If $x'<\XIDSTAR=\min(\XISSTAR)$ we observe that for the optimal intervention $\XIDSTAR$ we have $\SMINT{\ENDEPOCH}(\XIDSTAR)=0$ and $\XIDSTAR$ satisfies hypothesis of \PROREF{PN3:WorsenByDownShift}. Thus if we write $x'=\XIDSTAR+(x'-\XIDSTAR)$ we get $L(x')=\XIDSTAR-x'$ and $\ENDSTOCK{x'}=\ENDSTOCK{\XIDSTAR}$.

\medskip
    Finally, the last two steps proves that no other intervention outside $\XISSTAR$ can be optimal.
\end{itemize}
\end{myproof}

\begin{myproof}[Proof of \PROREF{PN7:XQfromXI}]
    Let us analyze the three cases.

    \medskip
    If $\XISTAR\in[q-Q,q]$, then $\XISTAR$ is also feasible for the more constrained capacitated problem. Equalities $\XQSTAR=\XISTAR$, $\LQSTAR=\LISTAR$ and $s^{(\ENDEPOCH)}(\XQSTAR)=s^{(\ENDEPOCH)}(\XISTAR)$ are thus obvious.
    Equality $s^{(\ENDEPOCH)}(\XISTAR)=s^{(\ENDEPOCH)}(0)$ comes from the following argument. If $\XISTAR=0$, there is nothing more to prove. Otherwise, the equality comes from \PROREF{PN6:MultipleOptimalInterventions} by setting $x'=0$ in both cases, i.e., $\LISTAR>0$ and $\LISTAR=0$, respectively. In the latter case use second implication when $\XIUSTAR=\XISTAR<0$, and third implication when $0<\XISTAR=\XIDSTAR$.

    \medskip
    If $q<\XISTAR$, then $0<\XISTAR=\XIDSTAR$ and any feasible intervention $x'\in[q-Q,q]$ for the capacitated problem is smaller than \XISTAR; according to \PROREF{PN6:MultipleOptimalInterventions} we get $L(x')=\LISTAR+|\XISTAR-x'|$ which is clearly minimum only for the maximum feasible value for $x'$. Thus, $\XQSTAR=q$ and $\LQSTAR=\LISTAR+|\XISTAR-\XQSTAR|$.
    Condition $s^{(\ENDEPOCH)}(\XQSTAR)=s^{(\ENDEPOCH)}(\XISTAR)=s^{(\ENDEPOCH)}(0)$ is
    directly derived from \PROREF{PN6:MultipleOptimalInterventions}) by setting $x'=\XQSTAR$ and $x'=0$ (use third implication in case 2).

    \medskip
    If $\XISTAR<q-Q$, then $\XIUSTAR=\XISTAR<0$ and any feasible intervention $x'\in[q-Q,q]$ for the capacitated problem is greater than \XISTAR. The argument follows the same line of the previous case.
\end{myproof}

\bigskip
\begin{myproof}[Proof of \TEOREF{TH2:alg:Vehicle-Intervention}]

\paragraph{Complexity.}
We assume information on the system state induced by the null intervention at all epochs $h\in[\STARTEPOCH,\ENDEPOCH]$, described in vectors $\SBF$, $\lmBF$ and $\lpBF$, has been already computed in $O(m)$ time according to equations \EQUREF{eq02:update_shat_int}-\EQUREF{eq06:update_s}.
The algorithm itself runs in $O(\ENDEPOCH-\STARTEPOCH)$, indeed every step has a cost $O(1)$ while the main loop executes $O(\ENDEPOCH-\STARTEPOCH)$ times.

\paragraph{Correctness.}
We show that algorithm \CiteAlgoVehicleIntervention\ initializes and sequentially updates the $4$-tuple of variables $\langle \XISTAR, \LISTAR, \DOWN, \UP\rangle$, so that, at the end of each loop iteration, the following property holds on interval $[\STARTEPOCH,h]$:
\begin{itemize}
    \item \XISTAR\ is the minimum modulus optimal intervention;
    \item \LISTAR\ is the optimal loss;
    \item $s^{(h)}(\XISTAR)=s^{(h)}(0)$;
    \item $\DOWN=\SMINT{h}(\XISTAR)$;
    \item $\UP=\SMAXT{h}(\XISTAR)$.
\end{itemize}
Let us observe that if the property holds, from \PROREF{PN5:OptimalityConditions} we derive quite straightforwardly that
$\LISTAR>0\Longrightarrow(\DOWN=0\wedge\UP=C)$, or viceversa that $(\DOWN>0\vee\UP<C)\Longrightarrow\LISTAR=0$.

\noindent
We proceed with the following steps.
\begin{enumerate}
    \item Property holds after the initialization phase (lines \LINEREF{step:firstupdateX} -- \LINEREF{step:firstupdateBOUNDS}) for $h=\STARTEPOCH$.
    Observe that the $4$-tuple is consistently initialized with $\langle \XISTAR=l^{(\STARTEPOCH)-}-l^{(\STARTEPOCH)+}, \LISTAR=0, \DOWN=s^{(\STARTEPOCH)}, \UP=s^{(\STARTEPOCH)}\rangle$.
    Indeed, if no loss is induced at epoch $\STARTEPOCH$ by the null intervention (i.e., $l^{(\STARTEPOCH)-}=l^{(\STARTEPOCH)+}=0$), then the minimum modulus optimal intervention is zero, the loss is null and the stock level at epoch \STARTEPOCH\ is unchanged.
    Otherwise let us assume  w.l.o.g. we have a stockout loss at epoch \STARTEPOCH\ (i.e., $l^{(\STARTEPOCH)-}>0$ and $l^{(\STARTEPOCH)+}=0$).
    In this case, we have $\XISTAR=l^{(\STARTEPOCH)-}$ which is the minimum positive intervention to recover the whole loss without changing the stock level $s^{(\STARTEPOCH)}(\XISTAR)=s^{(\STARTEPOCH)}(0)=0$ at epoch \STARTEPOCH;
    in fact, according to equations \EQUREF{eq02:update_shat_int} -- \EQUREF{eq06:update_s} a smaller intervention would leave a positive loss at epoch \STARTEPOCH, while a larger intervention would increase the stock level at epoch \STARTEPOCH.
    A similar argument holds in case $l^{(\STARTEPOCH)+}>0$.

    \item Assuming that for some $h>\STARTEPOCH$ the property holds at the beginning of a loop iteration on interval $[\STARTEPOCH,h-1]$, 
    we prove that at the end of the iteration the property holds on interval $[\STARTEPOCH,h]$.
    First, observe that
    since $s^{(h-1)}(\XISTAR)=s^{(h-1)}(0)$, we have for all $h'\geq h$
    \begin{align*}
    & s^{(h')}(\XISTAR)=s^{(h')}(0),\\
    & l^{(h')-}(\XISTAR)=l^{(h')-}(0),\\ 
    & l^{(h')+}(\XISTAR)=l^{(h')+}(0), 
    \end{align*}
    thus,
    when the iteration starts, $\LISTAR+l^{(h)+}+l^{(h)-}$ represents the loss induced on $[\STARTEPOCH,h]$ by \XISTAR, which is currently  the minimum modulus optimal solution on interval $[\STARTEPOCH,h-1]$ with optimal loss \LISTAR.
    Let us also observe that at least one of $\delta^{-}$ and $\delta^{+}$ must be null and proceed by cases:
    \begin{itemize}
        \item $\delta^{+}>0$ (and $\delta^{-}=0$).
        According to computation at line \LINEREF{step:recoversurplus},
        condition $\delta^{+}>0$ implies $0<\DOWN=\SMINT{h-1}(\XISTAR)$ and $l^{(h)+}(\XISTAR)>0$.

        Inequality $\DOWN>0$ implies $\LISTAR=0$. Indeed, $\SMINT{h-1}(\XISTAR)>0$ (minimum stock level in interval $[\STARTEPOCH,h-1]$) implies $h^{(-)}(\XISTAR)=\infty$ and both optimality conditions \EQUREF{PN5_10:OptimalityConditions_BC} and \EQUREF{PN5_11:OptimalityConditions_DE} in \PROREF{PN5:OptimalityConditions} are excluded.

        On the other hand, inequality $l^{(h)+}(\XISTAR)>0$ implies $s^{(h)}(\XISTAR)=s^{(h)}(0)=C$.

        In particular, from $\LISTAR=0$, we have $l^{(h')+}(\XISTAR)=l^{(h')-}(\XISTAR)=0\quad \forall h'\leq h-1$, which together with $l^{(h)+}(\XISTAR)>0$ implies $h=h^{(+)}(\XISTAR)$.
        Thus, $h$ and $\delta^{+}$ play the role of $\underline{h}$ and $\delta^*$ in \PROREF{PN1:ImproveByDownShift}.
        Accordingly, the intervention update at line \LINEREF{step:updateX} is an improving intervention on interval $[\STARTEPOCH,h]$ with a loss update given at line \LINEREF{step:updateL}.
        Moreover, \PROREF{PN1:ImproveByDownShift} guarantees that stock levels are downshifted up to epoch $h-1$, and unchanged starting from $h$ (i.e., $s^{(h')}(\XISTAR)=s^{(h')}(0)\quad\forall h'\geq h$) including $s^{(h)}(\XISTAR)=s^{(h)}(0)=C.$
        Then, variables $\DOWN$ and $\UP$
        are correctly updated to
        $\DOWN=\SMINT{h}(\XISTAR)$
        at line \LINEREF{step:updateD}
        and to
        $\UP=\SMAXT{h}(\XISTAR)=C$
        at line \LINEREF{step:updateU}.

        The optimality of the new value for \XISTAR\ in interval $[\STARTEPOCH,h]$ is argued as follows.
        If $\delta^+$ computed at line \LINEREF{step:recoversurplus} equals $l^{(h)+}$, we know that the updated value of \LISTAR\ is still zero and we cannot do any better.
        Otherwise, we have the first loss in $h=h^{(+)}(\XISTAR)$ while the new value of $\DOWN$ (as updated at line \LINEREF{step:updateD}) equals $0$, thus updated intervention \XISTAR\ matches optimality condition \EQUREF{PN5_10:OptimalityConditions_BC} in \PROREF{PN5:OptimalityConditions} with respect to interval $[\STARTEPOCH,h]$

        Finally, according to \PROREF{PN1:ImproveByDownShift} we know that any update to \XISTAR\ smaller than $\delta^+$ would imply a larger update to \LISTAR, which proves the suboptimality of any other intervention smaller in modulus.

        \item $\delta^{-}>0$ (and $\delta^{+}=0$).
        Discussion mirroring that of previous point.

        \item $\delta^{+}=\delta^{-}=0$.
        Intervention \XISTAR\ does not change and loss \LISTAR\ and bounds $\DOWN$ and $\UP$ are consistently updated at lines \LINEREF{step:updateL}, \LINEREF{step:updateD} and \LINEREF{step:updateU}, respectively, with $s^{(h)}(\XISTAR)=s^{(h)}(0)$ implied by $s^{(h-1)}(\XISTAR)=s^{(h-1)}(0)$.

        Observe that, by hypothesis, optimality conditions stated in \PROREF{PN5:OptimalityConditions} hold for \XISTAR\ on $[\STARTEPOCH,h-1]$.
        If either condition \EQUREF{PN5_10:OptimalityConditions_BC} or \EQUREF{PN5_11:OptimalityConditions_DE} is satisfied, then the same condition is still satisfied on interval $[\STARTEPOCH,h]$ (epoch of first loss does not change with respect to the larger interval and so does the stock level up to then).
        If only condition \EQUREF{PN5_09:OptimalityConditions_A} is satisfied then we have three possible cases.
        \begin{itemize}
            \item $l^{(h)+}=l^{(h)-}=0.$ Optimality condition \EQUREF{PN5_09:OptimalityConditions_A} still holds on interval $[\STARTEPOCH,h]$.

            \item $l^{(h)+}>0$, $\DOWN=0.$ On interval $[\STARTEPOCH,h]$ we have $\SMINT{h}=0$ and $h=h^{(+)}(\XISTAR)<h^{(-)}(\XISTAR)$, thus optimality condition \EQUREF{PN5_10:OptimalityConditions_BC} hold for \XISTAR\ on interval $[\STARTEPOCH,h]$.

            \item $l^{(h)-}>0$, $\UP=C.$ On interval $[\STARTEPOCH,h]$ we have $\SMAXT{h}=C$ and $h=h^{(-)}(\XISTAR)<h^{(+)}(\XISTAR)$, thus optimality condition \EQUREF{PN5_11:OptimalityConditions_DE} hold for \XISTAR\ on interval $[\STARTEPOCH,h]$.

        \end{itemize}
        
        So far we know that $\XISTAR$ is the minimum modulus optimal intervention on interval $[\STARTEPOCH,h-1]$, and an optimal intervention on $[\STARTEPOCH,h]$. 
        To show that $\XISTAR$ is also the minimum modulus optimal intervention on interval $[\STARTEPOCH,h]$ we use the following argument.
        If $\XISTAR=0$ the point is self evident. Otherwise, let us assume w.l.o.g. that $\XISTAR>0$ and keep the focus on interval $[\STARTEPOCH,h-1]$ for which we have $\XISTAR=\XIDSTAR>0$. Now consider an alternative intervention $x'<\XISTAR$. From \PROREF{PN6:MultipleOptimalInterventions} (use $h-1$ in place of $\ENDEPOCH$) we observe that $x'$ is suboptimal on interval $[\STARTEPOCH,h-1]$ and $s^{(h-1)}(x')=s^{(h-1)}(\XISTAR)$.
        Thus, according to equations \EQUREF{eq02:update_shat_int}-\EQUREF{eq06:update_s}, interventions $x'$ and $\XISTAR$ produce the same loss at epoch $h$. This proves the suboptimality of $x'$ on interval $[\STARTEPOCH,h]$.

    \end{itemize}

    \item Finally (lines \LINEREF{step:fix1} -- \LINEREF{step:fix3}), the algorithm computes $\XQSTAR$ and $\LQSTAR$ accordingly to \PROREF{PN7:XQfromXI}.
\end{enumerate}
\end{myproof}

\begin{myproof}[Proof of \PROREF{PN8:PALGO.S0}]

Let us proceed in order.
\begin{itemize}
    \item [Eq. \EQUREF{PN8_18:PALGO.SO.1}]
    We first observe that
    \begin{align}\label{PN8_proof26:EQ:Traslazione}
    \LI{s}{x}=\LI{s+\delta}{x-\delta}.
    \end{align}
    Indeed, for any intervention $x$ and $\delta$ we have
    \begin{align}
      \hat{s}^{(\STARTEPOCH)}(x-\delta\ |\ s+\delta)=(s+\delta)+(x-\delta)+d^{(\STARTEPOCH)}=s+x+d^{(\STARTEPOCH)}=\hat{s}^{(\STARTEPOCH)}(x\ |\ s)\label{PN8_proof27:EQ:STARTSTOCK},
    \end{align}
    that is, from epoch $\STARTEPOCH$ all values of vectors \SBF, \lmBF, \lpBF\ are the same. In particular we have the same losses all along the time interval and the same stock level at the final epoch \ENDEPOCH.

    Then we observe that for any $x$ we have an intervention $y=x-\delta$ such that $\LI{s}{x}=\LI{s+\delta}{y}$, thus $$\min_{x}\LI{s}{x}\geq\min_{x}\LI{s+\delta}{x}.$$
    Analogously, for any $x$ we have an intervention $y=x+\delta$ such that $\LI{s+\delta}{x}=\LI{s}{y}$, thus $$\min_{x}\LI{s}{x}\leq\min_{x}\LI{s+\delta}{x}$$  and clearly, $$\min_{x}\LI{s}{x}=\min_{x}\LI{s+\delta}{x}\DEFINE\LISTAR.$$

    Moreover, from \eqref{PN8_proof26:EQ:Traslazione} we get that $\LI{s}{x}=\LISTAR$ if and only if $\LI{s+\delta}{x-\delta}=\LISTAR$, thus  $\XIDSTAR(s+\delta)=\XIDSTAR(s)-\delta$ and $\XIUSTAR(s+\delta)=\XIUSTAR(s)-\delta$. This concludes the argument.

    \item [Eq. \EQUREF{PN8_19:PALGO.SO.4}]
    
    Let us consider $\LISTAR>0$ and fix $x$ and $s\in[0,C]$.
    Now let $\delta=-s$ and use equation \EQUREF{PN8_proof27:EQ:STARTSTOCK} to realize that
    $$s^{(\ENDEPOCH)}(x\ |\ s) = s^{(\ENDEPOCH)}(x+s\ |\ 0).$$
    Finally, from equations \EQUREF{PN6_12:MOI.1} of \PROREF{PN6:MultipleOptimalInterventions} we have
    $$s^{(\ENDEPOCH)}(x+s\ |\ 0)=s^{(\ENDEPOCH)}(0\ |\ 0),$$
    and obviously
    $$s^{(\ENDEPOCH)}(x\ |\ s)=s^{(\ENDEPOCH)}(0\ |\ 0).$$

    \item [Eq. \EQUREF{PN8_20:PALGO.SO.2}]
    Observe that the virtual loss induced at epoch $\STARTEPOCH$ by intervention $y=\XISTAR(s)-\delta\in[q,\XISTAR(s)]$ with initial stock level $s+\delta$ is, according to \EQUREF{PN8_proof27:EQ:STARTSTOCK},
    $$\hat{s}^{(\STARTEPOCH)}(y\ |\ s+\delta)=\hat{s}^{(\STARTEPOCH)}(\XISTAR(s)\ |\ s),$$
    that is $y$ produces, at every epoch in interval $[\STARTEPOCH,\ENDEPOCH]$  the same stock levels and losses.
    Thus, $y$ produces the same total loss \LISTAR\ and so is optimal for the initial stock level $s+\delta$.
    In the same way we can show that for any $y'\in[0,y)$ we get
    $$\hat{s}^{(\STARTEPOCH)}(y'\ |\ s+\delta)=\hat{s}^{(\STARTEPOCH)}(\XISTAR(s)-(y-y')\ |\ s)$$ and, because $\XISTAR(s)-(y-y')\in[0,\XISTAR(s))$, $y'$ produces a loss larger than $\LISTAR$, thus,
    $$y=\XISTAR(s)-\delta=\XISTAR(s+\delta).$$

    Equalities
    \begin{align*}
      & \XQSTAR(s+\delta) = \XQSTAR(s) = q  \qquad\qquad\text{and}\\
      & \LQSTAR(s+\delta) = \LISTAR+\XISTAR(s)-q-\delta
    \end{align*}
    come directly from condition $\XISTAR(s+\delta)=y>q$ and equations \EQUREF{PN7_15:PALGO.2.eq}, \EQUREF{PN7_16:PALGO.4.eq3} in \PROREF{PN7:XQfromXI}, respectively.

    Finally, observe that from $\XISTAR(s+\delta)=\XISTAR(s)-\delta$, and equation \EQUREF{PN8_proof27:EQ:STARTSTOCK} we have at epoch \STARTEPOCH
    $$\hat{s}^{(\STARTEPOCH)}(\XISTAR(s+\delta)\ |\ s+\delta)=\hat{s}^{(\STARTEPOCH)}(\XISTAR(s)-\delta\ |\ s+\delta)=\hat{s}^{(\STARTEPOCH)}(\XISTAR(s)\ |\ s),$$
    then we also have at epoch \ENDEPOCH\
    $$s^{(\ENDEPOCH)}(\XISTAR(s+\delta)\ |\ s+\delta)=s^{(\ENDEPOCH)}(\XISTAR(s)\ |\ s).$$
    Also consider that $\XQSTAR(s+\delta)=\XQSTAR(s)=q\leq \XISTAR(s+\delta)\leq \XISTAR(s)$.
    Thus, from either equations \EQUREF{PN6_12:MOI.1} or \EQUREF{PN6_13:MOI.2} of \PROREF{PN6:MultipleOptimalInterventions}, we derive
    \begin{align*}
        & s^{(\ENDEPOCH)}(\XISTAR(s+\delta)\ |\ s+\delta)= s^{(\ENDEPOCH)}(\XQSTAR(s+\delta)\ |\ s+\delta)\\
        & s^{(\ENDEPOCH)}(\XISTAR(s)\ |\ s)=s^{(\ENDEPOCH)}(\XQSTAR(s)\ |\ s)
    \end{align*}
    and so
    $$s^{(\ENDEPOCH)}(\XQSTAR(s+\delta)\ |\ s+\delta)=s^{(\ENDEPOCH)}(\XQSTAR(s)\ |\ s).
    $$

    \item [Eq. \EQUREF{PN8_21:PALGO.SO.3}]
    Discussion follows the same line as equalities \EQUREF{PN8_20:PALGO.SO.2}.
\end{itemize}
\end{myproof}

\subsection{General problem} 

\begin{myproof}[Proof of \PROREF{PN9:OttimoOgniIntervallo}]

We proceed by contradiction. Let $\XBF^*=\XVECTrow[x^*]{1}{i}$ be a vector of interventions and assume that for some $k\leq i$, intervention $x^*_k$ is not optimal for the \ONEINTPROBLEM\ capacitated problem on interval $I_k$ with initial stock level
$\SX[x^*]{k}$ at epoch $e_{k-1}$, i.e.,
$$\LIOTT[x^*]{k}>\underset{x\in[q_k-Q_k,q_k]}{\min}\ \LIXS[x]{k}{\SX[x^*]{k}}.$$
Now, let us focus on the \ONEINTPROBLEM\ problem on interval $I_k$.
Since $x^*_k$ is suboptimal for the capacitated problem, it must be either $x^*_k<\XQDSTAR$ or $\XQUSTAR<x^*_k$.
If $x^*_k<\XQDSTAR$, we must have $q-Q<\XQDSTAR\leq\XIDSTAR$ and we set $\widetilde{y}=\XQDSTAR$; otherwise we must have $\XIUSTAR\leq\XQUSTAR<q$ and we set $\widetilde{y}=\XQUSTAR$. In any case $\widetilde{y}$ is optimal for the capacitated \ONEINTPROBLEM\ problem; moreover we have
$$s^{(e_{k+1}-1)}([x^*_1,\ldots,x^*_{k-1},\widetilde{y}])=\SXPIU[x^*]{k}$$
where the equality is guaranteed by \PROREF{PN6:MultipleOptimalInterventions} (both $x^*_k$ and $\widetilde{y}$ induce the same stock level as the optimal intervention $\XIDSTAR$ or $\XIUSTAR$).

Thus, for vector of interventions $\YBF=\XVECTrow[y]{1}{i}$ with $y_j=x^*_j$ for $j\neq k$ and $y_k=\widetilde{y}$ we have from equation \EQUREF{eq22:EQ:sumLIOTT}
\begin{align*}\begin{array}{rcl}
\LVECTBKW[y]{i} & =& \sum_{h=1}^{i}\LIOTT[y]{h}\\
                & <& \sum_{h=1}^{i}\LIOTT[x^*]{h}=\LVECTBKW[x^*]{i}
\end{array}
\end{align*}
which contradicts optimality of $\XBF^*$ for problem \PROBLEMA{i}.
\end{myproof}

\begin{myproof}[Proof of \PROREF{PN10:ottimoIntervalli}]
\newcommand{\maxLXY}[1][h']{\ensuremath{\max\left(\LIOTT[x^*]{#1},\LIOTT[y^*]{#1}\right)}}

Property is guaranteed by the two following claims for all $h\leq i$
\begin{enumerate}[A.]
    \item\label{CLAIM:A} $\SX[x^*]{h}=\SX[y^*]{h} \Longrightarrow \LIOTT[x^*]{h} = \LIOTT[y^*]{h}$;
    \item\label{CLAIM:B} $\SX[x^*]{h}\neq\SX[y^*]{h} \Longrightarrow \maxLXY[h]=\LISTAR_{I_{h}}$.
\end{enumerate}
Claim \ref{CLAIM:A} is guaranteed by  \PROREF{PN9:OttimoOgniIntervallo}. 
With respect to claim \ref{CLAIM:B}, observe that, according to equation \EQUREF{PN7_16:PALGO.4.eq3} of \PROREF{PN7:XQfromXI}, $\LISTAR_{I_h}$ is a lower bound for both $\LIOTT[x^*]{h}$ and $\LIOTT[y^*]{h}$, thus the equality also implies $\LIOTT[x^*]{h}=\LIOTT[y^*]{h}$.

\medskip
We prove claim \ref{CLAIM:B} by absurdity. 
We assume that an index $h\leq i$ exists such that 
$$\SX[x^*]{h}\neq\SX[y^*]{h} \text{ and } \maxLXY[h]>\LISTAR_{I_h},$$
and get to some contradiction.

Let $h'$ be the minimum of such indices, and let $l$ be the last index before $h'$ such that $\SX[x^*]{l}=\SX[y^*]{l}$.
Note that such an index $l$ exists as the equality holds at least for $l=1$ (possibly $l=h'-1$), and observe that we have the following Facts:
\medskip
\begin{enumerate}[(a)]
    
    \item\label{FACT:A} $\SX[x^*]{l}=\SX[y^*]{l}$ and $\SX[x^*]{k}\neq\SX[y^*]{k}$ for $k=l+1,\ldots,h'$; by definition of $h'$ and $l$.
    
    \item\label{FACT:B} $\LIOTT[x^*]{k}=\LIOTT[y^*]{k}=\LISTAR_{I_k}$ for $k=l+1\ldots,h'-1$; by definition of $h'$ and $l$.
    
    \item\label{FACT:D} $\LIOTT[x^*]{l}=\LIOTT[y^*]{l}=\LISTAR_{I_l}=0$; 
    this follows from \PROREF{PN9:OttimoOgniIntervallo}, 
    which guarantees that $x^*_l$ and
    $y^*_l$ are both optimal for the 
    same \ONEINTPROBLEM\ capacitated problem, 
    from Fact \FACTREF{FACT:A} which implies $x^*_{l}\neq y^*_{l}$, 
    and from equation \EQUREF{PN7_14:PALGO.2.eqXR} of \PROREF{PN7:XQfromXI} showing that two distinct optimal intervention for the same \ONEINTPROBLEM\ capacitated problem are allowed only if the optimal loss is zero.
    
    \item\label{FACT:E} $\LISTAR_{I_k}=0$ for $k=l,\ldots,h'-1$; this follows from equation \FACTREF{PN8_19:PALGO.SO.4} of \PROREF{PN8:PALGO.S0}, which guarantees that $\LISTAR_{I_k}>0$ implies $\SX[x^*]{k+1}=\SX[y^*]{k+1}$ in contradiction with Fact \FACTREF{FACT:A}.
    
    \item\label{FACT:F} $\LIOTT[x^*]{k}=\LIOTT[y^*]{k}=0$ for $k=l,\ldots,h'-1$; this is a straight consequence of 
    \PROREF{PN9:OttimoOgniIntervallo}, which guarantees that $x^*_k$ and $y^*_k$ are both optimal for the corresponding \ONEINTPROBLEM\ capacitated problem on intervals $I_k$, and from
    Facts \FACTREF{FACT:B}, \FACTREF{FACT:D} and \FACTREF{FACT:E}.
    
    \item\label{FACT:G}
                $\SX[y^*]{k}-\SX[x^*]{k}=\sum_{j=l}^{k-1}(y^*_j-x^*_j)$ for $k=l+1,\ldots,h'$;
            this follows from Fact \FACTREF{FACT:F} and equations \EQUREF{eq02:update_shat_int} -- \EQUREF{eq06:update_s}. 
\end{enumerate}

\bigskip

Now, let us assume w.l.o.g. that $\SX[y^*]{h'}>\SX[x^*]{h'}$ and discuss the possible cases.
\begin{itemize}
    \item 
$\maxLXY=\LIOTT[x^*]{h'}>\LISTAR_{I_{h'}}$. Then, $x^*_{h'}\in\XQSSTAR_{I_{h'}}(\SX[x^*]{h'})$, but $x^*_{h'}\notin\XISSTAR_{I_{h'}}(\SX[x^*]{h'})$ and, according to
    equation \EQUREF{PN7_15:PALGO.2.eq}
    of \PROREF{PN7:XQfromXI}, we have two possible cases:
\begin{itemize}
    \item 
    Case $x^*_{h'}=q_{h'}<\XIDSTAR_{I_{h'}}(\SX[x^*]{h'})$. According to Fact \EQUREF{FACT:G},
        inequality $\SX[y^*]{h'}>\SX[x^*]{h'}$ implies that at least one index $k<h'$ exists such that $y^*_k>x^*_k$,
        and by defining $\hat{k}$ as the last of such indices ($l\leq\hat{k}<h'$), 
        we have that $\SX[y^*]{k}-\SX[x^*]{k}$ is monotonically non-increasing with respect to $k$ starting from $\hat{k}+1$ up to $h'$.
        
        Now consider the vector of interventions $\ZBF$ with 
        $z_k=x^*_k$ for $k\neq \hat{k}$ and 
        $z_{\hat{k}}=x^*_{\hat{k}}+\delta$ for 
        $\delta=1\leq \min(\SX[y^*]{h'}-\SX[x^*]{h'},\XIDSTAR_{I_{h'}}(\SX[x^*]{h'})-q_{h'})$; 
        observe that from equations \EQUREF{eq02:update_shat_int} -- \EQUREF{eq06:update_s} we easily derive the following Facts:
        \begin{enumerate}[(i)]
            \item\label{FACT:I} $\SX[z]{k}=\SX[x^*]{k}$ and $\LIOTT[z]{k}=\LIOTT[x^*]{k}$ for $k=1,\ldots,\hat{k}$;
            \item\label{FACT:II} $\SX[z]{k}=\SX[x^*]{k}+\delta\in[\SX[x]{k},\SX[y]{k}]$ for $k=\hat{k}+1,\ldots,h'$;
            \item\label{FACT:III} 
            $\LIOTT[z]{k}=\LIOTT[x^*]{k}=0$  for $k=\hat{k},\ldots,h'-1$;
        \end{enumerate}

\medskip
As for $\LIOTT[z]{h'}$ and $\LIOTT[x^*]{h'}$ we observe that from Fact \EQUREF{FACT:II} we have
$$\SX[z]{h'}=\SX[x^*]{h'}+\delta,$$
and according to equations \EQUREF{PN8_20:PALGO.SO.2} of \PROREF{PN8:PALGO.S0}, we get
\begin{align*}
    \begin{array}[l]{lll}
    z_{h'}=q_{h'}=\XQSTAR_{I_{h'}}(\SX[z]{h'}),\\
    \begin{array}[l]{lll}
    \LIOTT[z]{h'} & =\LQSTAR_{I_{h'}}(\SX[z]{h'})\\
                       & =\LQSTAR_{I_{h'}}(\SX[x^*]{h'})-\delta\\
                       & =\LIOTT[x^*]{h'}-\delta\\
                       & <\LIOTT[x^*]{h'},
    \end{array}\\
    \SX[z]{h'+1}=\SX[x^*]{h'+1}.
    \end{array}
\end{align*}
The first equality shows that $z_{h'}=q_{h'}$ is the optimal intervention for the \ONEINTPROBLEM\ capacitated problem on interval $I_{h'}$ given the initial stock level $\SX[z]{h'}$. 
The second equality implies that vector of interventions \ZBF\ attains a smaller loss in interval $I_{h'}$ with respect to $\XBF^*$.
The third equality guarantees that starting from epoch
$e_{h'+1}$ onward the stock level induced by vector \ZBF\ is unchanged with respect to $\XBF^*$, and so are the losses.
In conclusion, we have $\LVECTBKW[z]{i}<\LVECTBKW[x^*]{i}$, in contradiction with the optimality of $\XBF^*$ for \PROBLEMA{i}.

\item Case $x^*_{h'}=q_{h'}-Q_{h'}>\XIUSTAR_{I_{h'}}(\SX[x^*]{h'})$. 
We first observe that from equations \EQUREF{PN8_18:PALGO.SO.1} of \PROREF{PN8:PALGO.S0} we obtain
\begin{align*}
\XIUSTAR_{I_{h'}}(\SX[y^*]{h'}) & =\XIUSTAR_{I_{h'}}(\SX[x^*]{h'})-\left(\SX[y^*]{h'}-\SX[x^*]{h'}\right)\\
                                & <\XIUSTAR_{I_{h'}}(\SX[x^*]{h'})\\
                                & <x^*_{h'}=q_{h-Q_{h'}}\leq0    
\end{align*}

and derive that $y^*_{h'}=\XQUSTAR_{I_{h'}}(\SX[y^*]{h'})=q_{h'}-Q_{h'}=x^*_{h'}$.

Moreover, note that $\XIUSTAR_{I_{h'}}(\SX[y^*]{h'})-(q_{h'}-Q_{h'})<\SX[x^*]{h'}-\SX[y^*]{h'}<0$.
Thus, setting $\delta=\SX[x^*]{h'}-\SX[y^*]{h'}$ and $s=\SX[y^*]{h'}$ we can obtain 
from equations \EQUREF{PN8_21:PALGO.SO.3} of \PROREF{PN8:PALGO.S0} the following chain of inequalities
        \begin{align*}
            \LIOTT[x^*]{h'} & =\LQSTAR_{I_{h'}}(\SX[x^*]{h'})\\
                           & =\LQSTAR_{I_{h'}}(\SX[y^*]{h'}+\SX[x^*]{h'}-\SX[y^*]{h'})\\
                           & =\LQSTAR_{I_{h'}}(\SX[y^*]{h'})+\left(\SX[x^*]{h'}-\SX[y^*]{h'}\right)\\
                           & <\LQSTAR_{I_{h'}}(\SX[y^*]{h'}) =\LIOTT[y^*]{h'},
        \end{align*}
in contradiction with $\LIOTT[x^*]{h'}\geq\LIOTT[y^*]{h'}$.

\end{itemize}
\item 
$\maxLXY=\LIOTT[y^*]{h'}>\LISTAR_{I_{h'}}$. Discussion mirroring that of previous point.
\end{itemize}
This concludes the argument.
\end{myproof}

\begin{myproof}[Proof of \PROREF{PN11:OttimoTotaleParziale}]

We proceed by induction on the dimension $i$ of the optimal vector \XBF\ for \PROBLEMA{i}.
The inductive hypothesis is that an optimal vector of interventions for \PROBLEMA{i} is also optimal for \PROBLEMA{h} for all $1\leq h< i$.

\paragraph{Base case.}
Let $i=2$ and $\XBF=[x_1,x_2]$ be an optimal vector of interventions for \PROBLEMA{i}. We show that \XBF\ satisfies the inductive hypothesis.
    \PROREF{PN9:OttimoOgniIntervallo} guarantees that intervention $x_1$ is optimal for the \ONEINTPROBLEM\ problem on interval $I_1$ with initial stock level $\SZERO$; but as $I_1=J_1$ vector $[x_1]$ is also optimal for problem \PROBLEMA{1}.

\paragraph{Inductive step.}
Let us assume that the inductive hypothesis holds for some $i\geq2$ and show that also holds for $i+1$.

    We proceed by contradiction.
    Let $\XBF=\XVECTrow[x]{1}{i+1}$ be an optimal vector of interventions for \PROBLEMA{i+1} which, for some $h\leq i$, is not optimal for \PROBLEMA{h}, and let $k$ be the minimum of such indices. 

    Let \YBF\ be any optimal vector of interventions for \PROBLEMA{i} and, according to the inductive hypothesis, also for \PROBLEMA{h}, $h=1,\ldots,i-1$.

    By construction of $k$, \XBF\ is optimal for \PROBLEMA{k-1} and sub-optimal for \PROBLEMA{k}.
    On the other side, by inductive assumption, \YBF\ is optimal for both \PROBLEMA{k-1} and \PROBLEMA{k}.
    Thus, from equation \EQUREF{eq23:EQ:LBCKWLIOTT} it is immediate to derive
    \begin{align*}
        \LIOTT[x]{k}>\LIOTT[y]{k},
    \end{align*}
    which, on the basis of \PROREF{PN9:OttimoOgniIntervallo}, excludes the case $\SX[x]{k}=\SX[y]{k}$.
    We conclude the discussion observing that the latter two conditions contradict the statement of \PROREF{PN10:ottimoIntervalli}.
\end{myproof}

\bigskip
\begin{myproof}[Proof of \PROREF{PN12:neutralStockOptimum}]
We consider only the case $i=\NUMV$. Discussion for $i<\NUMV$ follows the same line, just on a shorter time horizon.
If $\LVECTBKWO{\NUMV}=0$ then $\OBF$ is optimal and we have nothing to show.
Otherwise, let $\overline{h}$ be the last epoch when the null vector induces a loss
(i.e., $\overline{h}=\max\{h\ |\ L^{(h)}(\OBF)>0\}$) and let $\YBF$ be an optimal vector of interventions.

If $s^{(\overline{h})}(\YBF)=s^{(\overline{h})}(\OBF)$ then we obtain the optimal vector of interventions $\XBF^*$ by setting $x_j^*=y_j$ for $e_j\leq \overline{h}$ and $x_j^*=0$ for $e_j>\overline{h}$.
Observe that $\XBF^*$ produces the same effect (stock level and loss) as $\YBF$ up to epoch $\overline{h}$, and the same effect as $\OBF$ (with null loss) from $\overline{h}$ onward, which implies
$\LVECTBKW[x^*]{\NUMV}\leq\LVECTBKW[y]{\NUMV}$
and $s^{(m)}(\XBF^*)=s^{(m)}(\OBF)$.

\smallskip
To examine the case $s^{(\overline{h})}(\YBF)\neq s^{(\overline{h})}(\OBF)$, let us assume w.l.o.g. that we have a surplus loss at epoch $\overline{h}$ (i.e., $l^{(\overline{h})+}(\OBF)>0$) and necessarily $s^{(\overline{h})}(\YBF)<s^{(\overline{h})}(\OBF)=C$.

Let $e_{\overline{k}}$ be the last intervention epoch before $\overline{h}$ and
$e_k$ the first intervention epoch such that $s^{(\overline{h})}([y_1,\ldots,y_k,0_{k+1},\ldots,0_{\overline{k}}])<C$.
Clearly we have $e_k\leq e_{\overline{k}}\leq\overline{h}$, $s^{(\overline{h})}([y_1,\ldots,y_l,0_{l+1},\ldots,0_{\overline{k}}])=C$ for all $l<k$ and $y_k\neq0$.

\smallskip
Focusing on the \ONEINTPROBLEM\ problem on interval $[e_k,\overline{h}]$ with stock level \SX[y]{k} at epoch $e_k-1$ we observe that the stock level induced at epoch $\overline{h}$ by the null intervention at epoch $e_k$ is $C$ while the stock level induced by intervention $y_k$ is smaller than $C$.
According to \PROREF{PN6:MultipleOptimalInterventions} we can have two interventions producing different stock levels at the end of the interval only when the optimal loss for the uncapacitated problem is zero.
Moreover, according to equation \EQUREF{PN7_17:PALGO.4.eq4} of \PROREF{PN7:XQfromXI} minimum modulus optimal interventions \XISTAR\ and \XQSTAR\ for the uncapacitated and capacitated problems induce at the end of the interval the same stock level as the null intervention.
Summing up, 
we have
$L_{[e_k,\overline{h}]}\left([\XQSTAR,0_{k+1},\ldots,0_{\overline{k}}]\ \left|\ \SX[y]{k}\right.\right)=0$ with stock level $C$ at epoch $\overline{h}$.

Finally, observe that starting with a stock level equal to $C$ at epoch $\overline{h}$ the sequence of null interventions at epochs $e_i$ with $i>\overline{k}$ produces exactly the same effect, in terms of losses and stock levels, on the interval $[\overline{h},\NUMV]$; in particular, by construction of $\overline{h}$ the total loss induced in this interval is zero.

\medskip
Now, let us define a vector of interventions $\XBF$ as follows:
$x_j=y_j$ for $j<k$,
$x_k=\XQSTAR$,
and $x_j=0$ for $j>k$.
We are ready to show that $\XBF$ is optimal.

Indeed,
\begin{align*}\begin{array}{rcl}
\LVECTBKW[x]{\NUMV} & = & L_{[e_1,e_k-1]}\left([y_1,\ldots,y_{k-1}]\ \left|\ \SZERO\right.\right) +\\
& & L_{[e_k,\overline{h}]}\left([\XQSTAR,0_{k+1},\ldots,0_{\overline{k}}]\ \left|\ \SX[y]{k}\right.\right) +\\
& & L_{[\overline{h}+1,\NUMV]}\left(\OBF\ |\ C\right)\\
& = & L_{[e_1,e_k-1]}\left([y_1,\ldots,y_{k-1}]\ \left|\ \SZERO\right.\right)\\
& \leq & \LVECTBKW[y]{\NUMV}.
\end{array}
\end{align*}

Finally, condition $s^{(m)}(\XBF)=s^{(m)}(\OBF)$ comes by observing that by construction $s^{(\overline{h})}(\XBF)=s^{(\overline{h})}(\OBF)$ and starting from $\overline{h}$ all interventions in $\XBF$ are zero.
\end{myproof}

\bigskip
\begin{myproof}[Proof of \PROREF{PN13:augmentation}]
We first observe that for a given $i<\NUMV$, we have by construction
\begin{equation}\label{PN13_proof28:eq:19}
\LVECTBKW[x^*]{i}=\LVECTBKW[x^*]{i-1}+\LIOTT[x^*]{i}
\end{equation}
where $\XBF^*=\XVECTrow[x^*]{1}{i}$ is an optimal vector for the
problem \PROBLEMA{i}, and
\begin{equation}\label{PN13_proof29:eq:20}
\LVECTBKWXAUGMN{\bar{x}}{i}=\LVECTBKW[\bar{x}]{i-1}+\LIOTTAUGMN[\bar{x}]{i}{\SX[\bar{x}]{i}}
\end{equation}
where $\XSTARAUGMN=\XVECTrow[\bar{x}]{1}{i}$ is the optimal vector of interventions for the augmented problem $\PROBLEMAAUG{i}$ and the last term is the loss in the augmented interval $\AUGMINTI$ depending on the intervention $\bar{x}_{i}$ implemented at epoch $e_i$ and on the stock level \SX[\bar{x}]{i} determined at epoch $e_i-1$ by all previous interventions.

Here we assume that $\XBF^*$ is an optimal vector that satisfies condition stated in \PROREF{PN12:neutralStockOptimum} (i.e., $s^{(f)}(\XBF^*)=s^{(f)}(\OBF)$ for $f=e_{i+1}-1$) so that
\begin{equation}\label{PN13_proof30:eq:21}
\LVECTBKWXAUGMN{x^*}{i} = \LVECTBKW[x^*]{i}+|\delta|.
\end{equation}
Indeed, as $s^{(f)}(\XBF^*)=s^{(f)}(\OBF)$ for $f=e_{i+1}-1$, the total loss obtained in the augmented interval $\AUGMINTJ$ is given by the loss obtained in interval $J_i$ plus the loss obtained at epoch $e_{i+1}$ which is $|\delta|$ by construction.

Now, \PROREF{PN11:OttimoTotaleParziale} guarantees that partial vectors $\XVECTrow[\bar{x}]{1}{i-1}$ and $\XVECTrow[x^*]{1}{i-1}$ are optimal on interval $J_{i-1}$ and we have
\begin{equation}\label{PN13_proof31:eq:24}
\LVECTBKW[\bar{x}]{i-1}=\LVECTBKW[x^*]{i-1}.
\end{equation}

Thus, to prove optimality of $\XSTARAUGMN$ for interval $J_i$ it suffices to show that
$$\LIOTT[\bar{x}]{i}\leq\LIOTT[x^*]{i}.$$

\medskip\noindent
We prove this latter inequality by contradiction.
If
\begin{equation}\label{PN13_proof32:eq:23}
\LIOTT[\bar{x}]{i}>\LIOTT[x^*]{i},
\end{equation}
then there must be at least one epoch $h\in I_i$ with a positive loss, which according to equations
 \EQUREF{PN6_12:MOI.1} of \PROREF{PN6:MultipleOptimalInterventions}, does not change the stock and loss (w.r.t. $\OBF$) at later epochs and thus
we have
\begin{equation}\label{PN13_proof33:eq:22}
\LIOTTAUGMN[\bar{x}]{i}{\SX[\bar{x}]{i}}
= \LIOTT[\bar{x}]{i} +|\delta|.
\end{equation}

Using equations \EQUREF{PN13_proof28:eq:19} -- \EQUREF{PN13_proof33:eq:22} we can derive the following chain of inequalities

    \begin{tabular}{llr}
    \multicolumn{2}{l}{$\LVECTBKWXAUGMN{\bar{x}}{i} =$} & [by eq. \EQUREF{PN13_proof29:eq:20}]\\
          & $=\LVECTBKW[\bar{x}]{i-1}+\LIOTTAUGMN[\bar{x}]{i}{\SX[\bar{x}]{i}} =$ & [by eq. \EQUREF{PN13_proof33:eq:22}]\\
          & $= \LVECTBKW[\bar{x}]{i-1} + \LIOTT[\bar{x}]{i} +|\delta| =$ & [by eq. \EQUREF{PN13_proof31:eq:24}]\\
          & $= \LVECTBKW[x^*]{i-1} + \LIOTT[\bar{x}]{i} +|\delta| > $ & [by eq. \EQUREF{PN13_proof32:eq:23}]\\
          & $> \LVECTBKW[x^*]{i-1}+\LIOTT[x^*]{i}+|\delta| = $ & [by eq. \EQUREF{PN13_proof28:eq:19}]\\
          & $= \LVECTBKW[x^*]{i}+|\delta| = $ & [by eq. \EQUREF{PN13_proof30:eq:21}]\\
          & $=\LVECTBKWXAUGMN{x^*}{i}$
    \end{tabular}

proving that $\XSTARAUGMN$ cannot be an optimal vector of interventions for the augmented problem, in contradiction with the hypothesis.
\end{myproof}

\begin{myproof}[Proof of \PROREF{PN14:MERGE}]
We analyze the augmented problem \PROBLEMAAUG{i-1} proceeding by cases with respect to the sign of $\delta=\XISTAR(S_i)-\XQSTAR(S_i)$.
For sake of clarity the fictitious demand at epoch $e_i$ used to build \PROBLEMAAUG{i-1}\ is denoted by $\bar{d}_i$, not to be confused with $d^{(e_i)}$ which is the demand at epoch $e_i$ for problem \PROBLEMA{i}. 
According to equations \EQUREF{eq02:update_shat_int} -- \EQUREF{eq06:update_s}, 
we indicate with $\bar{s}^{(e_i)}(\YBF)=\max(0,\min(C,\SX[y]{i}+\bar{d}_i))$ the stock level induced by a vector of decisions \YBF\ at epoch $e_i$ in the augmented problem.

\paragraph{Case $\delta>0$.} 
When $\delta>0$ we know the following Facts:
\begin{enumerate}[(a)]

    \item\label{FACT:1}$\XISTAR(S_i)>\XQSTAR(S_i)=q_i$ (from \PROREF{PN7:XQfromXI}).
    
    \item\label{FACT:2}$\bar{d}_i=-\SX[0]{i}-\delta=-S_i-\delta$ and $\bar{s}^{(e_i)}(\OBF)=0$ (by construction of the augmented problem).

    \item\label{FACT:0}
    An optimal vector of interventions $\XSTARAUGMN$ for the augmented problem \PROBLEMAAUG{i-1}\ exists such that $\bar{s}^{(e_i)}(\bar{\XBF})=\bar{s}^{(e_i)}(\OBF)=0$ (from \PROREF{PN12:neutralStockOptimum} and Fact \eqref{FACT:2}).
    
    \item\label{FACT:5} \XSTARAUGMN\ is an optimal vector of interventions for problem \PROBLEMA{i-1} (from \PROREF{PN13:augmentation}).

    \item\label{FACT:6}
    Any optimal vector of interventions \YBF\ for \PROBLEMA{i-1} satisfies the equation
        \begin{align}
        \LVECTBKWXAUGMN{y}{i-1} &= \LSTAR[i-1]+\max(0,\delta-(\SX[y]{i}-S_i)).\label{PN14_proof34:FACT:6b}
        \end{align}
        Equation \EQUREF{PN14_proof34:FACT:6b} holds by construction of augmented problem, equation \eqref{eq05:update_lstockout} and Fact \EQUREF{FACT:2} applied to virtual stock of the augmented problem at epoch $e_i$: $\SX[y]{i}+\bar{d}_i=\SX[y]{i}-S_i-\delta$.

    \item\label{FACT:9}
    $\Delta\DEFINE(\LSTAR[i-1]+\delta) -  \LSTARAUGMN[i-1] \in [0,\delta]$
    (directly from $\LSTARAUGMN[i-1]\in \left[\LSTAR[i-1],\LSTAR[i-1]+\delta\right]$ where 
    the lower bound comes from equation \EQUREF{PN14_proof34:FACT:6b}; and
    the upper bound derives from equation \EQUREF{PN14_proof34:FACT:6b} applied to the vector \YBF\ that, 
    according to \PROREF{PN12:neutralStockOptimum}, is optimal for \PROBLEMA{i-1} with $\SX[y]{i}=S_i$).

    \item\label{FACT:10}
    $\SX[\bar{x}]{i} = S_i+\Delta$.
    (according to: Fact \eqref{FACT:9}, Fact \eqref{FACT:5}, equation \EQUREF{PN14_proof34:FACT:6b} applied to $\bar{\XBF}$, and Fact \EQUREF{FACT:0} that guarantees $\SX[\bar{x}]{i}+\bar{d}_i\leq0$; we have $\Delta =\LSTAR[i-1]+\delta -  \LSTARAUGMN[i-1]= \delta - \max(0,\delta-(\SX[\bar{x}]{i}-S_i))= \SX[\bar{x}]{i}-S_i$).
        
\end{enumerate}

We also observe that following Fact \EQUREF{FACT:10}, and according to Fact \EQUREF{FACT:1}, equation \EQUREF{PN8_20:PALGO.SO.2} of \PROREF{PN8:PALGO.S0} and equation \EQUREF{PN7_16:PALGO.4.eq3} of \PROREF{PN7:XQfromXI} we have 
\begin{align}\label{PN14_proof35}\begin{array}{lll}
&\XQSTAR(S_i+\Delta)&=\XQSTAR(S_i)=q_i,\\
&\XISTAR(S_i+\Delta)&=\XISTAR(S_i)-\Delta,\\
&\LQSTAR_{I_i}(S_i+\Delta)&=\LISTAR_{I_i}+(\XISTAR(S_i)-\Delta)-\XQSTAR(S_i)=\LISTAR_{I_i}+\delta-\Delta.
\end{array}\end{align}
Now, let us consider $\XBF^*=[\bar{x}_1,\ldots,\bar{x}_{i-1},q_i]$ and an optimal vector of interventions \YBF\ for \PROBLEMA{i}. 
Observe that, by construction, $\XBF^*$ is optimal for \PROBLEMA{i-1} and $\LVECTBKW[x^*]{i}=\LSTAR[i-1]+\LIOTT[x^*]{i}$ with $\LIOTT[x^*]{i}=\LQSTAR_{I_i}(S_i+\Delta)$.
On the other hand \YBF\ is optimal also for \PROBLEMA{i-1} and $\LVECTBKW[y]{i}=\LSTAR[i-1]+\LIOTT[y]{i}$.

Thus, to prove the optimality of $\XBF^*$ for \PROBLEMA{i} it suffices to show that $\LIOTT[x^*]{i}=\LIOTT[y]{i}$.
Let us proceed by cases.
\begin{itemize}
    \item If $\Delta=\delta$, then from \EQUREF{PN14_proof35} we have $\LIOTT[x^*]{i}=\LISTAR_{I_i}$ while, according to \PROREF{PN7:XQfromXI} and \PROREF{PN8:PALGO.S0}, we get $\LISTAR_{I_i}\leq\LIOTT[y]{i}$ where a strict inequality would be in contradiction with the optimality of \YBF.
    
    \item If $\Delta<\delta$, let us indicate with $\SX[y]{i}=S_i+\gamma$ the stock level induced by \YBF\ at epoch $e_i-1$.

    \begin{itemize}
        \item If $\gamma=\Delta$, then $\LIOTT[x^*]{i}= \LQSTAR_{I_i}(S_i+\Delta)$ by construction and $\LQSTAR_{I_i}(S_i+\Delta)=\LIOTT[y]{i}$ from \PROREF{PN9:OttimoOgniIntervallo}. 
        \item If $\gamma<\Delta$, then from equation \EQUREF{PN8_20:PALGO.SO.2} of \PROREF{PN8:PALGO.S0} we get $\LIOTT[y]{i}>\LIOTT[\bar{x}]{i}$ in contradiction with optimality of $\YBF$ for \PROBLEMA{i}.
        \item If $\gamma>\Delta$, observe that, according to equation \EQUREF{PN14_proof34:FACT:6b} we have $\LVECTBKWXAUGMN{x^*}{i-1} - \LVECTBKWXAUGMN{y}{i-1} = \delta-\Delta - \max(0,\delta-\gamma)>0$ which is in contradiction with the optimality of $\bar{\XBF}$ for the augmented problem \PROBLEMAAUG{i-1}.
    \end{itemize}
\end{itemize}

Finally, equality $\LSTAR[i]=\LSTARAUGMN[i-1]+\LISTAR_{I_i}$ straightly derives from equalities
$\LSTAR[i]=\LVECTBKW[x^*]{i}=\LSTAR[i-1]+\LQSTAR_{I_i}(S_i+\Delta)=\LSTAR[i-1]+\LISTAR_{I_i}+\delta-\Delta$
and the definition of $\Delta$.

\paragraph{Case $\delta<0$.} When $\delta<0$ we we can carry out a discussion mirroring the previous one.

\paragraph{Case $\delta=0$.} 
When $\delta=0$ we have by construction $\bar{d}_i=0$ and 
$\XISTAR(S_i)=\XQSTAR(S_i)$, hence $\LQSTAR_{I_i}(S_i)=\LISTAR_{I_i}$.
Let $\bar{\XBF}=\XVECTrow[\bar{x}]{1}{i-1}$ be the optimal vector of interventions for augmented problem \PROBLEMAAUG{i-1} such that $\bar{s}^{(e_i)}(\bar{\XBF})=\bar{s}^{(e_i)}(\OBF)$ as guaranteed by \PROREF{PN12:neutralStockOptimum}. Observe that not only is vector $\bar{\XBF}$ optimal for \PROBLEMA{i-1} as  guaranteed by \PROREF{PN13:augmentation}, but it also satisfies 
\begin{align}
   & \SX[\bar{x}]{i}=S_i,\label{PN14_36:augmntEqualNormal}\\
   & \LVECTBKWXAUGMN{\bar{x}}{i-1} =\LVECTBKW[\bar{x}]{i-1}= \LSTAR[i-1];\label{PN14_37:augmntEqualNormal}
\end{align}
equality \EQUREF{PN14_36:augmntEqualNormal} straightly derives from 
conditions $\bar{s}^{(e_i)}(\bar{\XBF})=\bar{s}^{(e_i)}(\OBF)$ and
$\bar{d}_i=0$, and from equations \EQUREF{eq02:update_shat_int}-\EQUREF{eq06:update_s} applied to the augmented problem;
equality \EQUREF{PN14_37:augmntEqualNormal} comes from $\bar{d}_i=0$ and optimality of $\bar{\XBF}$ for \PROBLEMA{i-1}.

Morever, 
vector $\XBF^*=[\bar{x}_1,\ldots,\bar{x}_{i-1},\XQSTAR(S_i)]$ produces on interval $J_i$ a loss
\begin{align*}
\LVECTBKW[x^*]{i} 
    &=\LSTAR[i-1] + \LIOTT[x^*]{i}\\
    &=\LSTAR[i-1] + \LQSTAR_{I_i}(S_i)\text{ (from $\SX[x^*]{i}=\SX[\bar{x}]{i}$, equation \EQUREF{PN14_36:augmntEqualNormal}, and $x^*_i=\XQSTAR(S_i)$)}\\
    &=\LSTAR[i-1] + \LISTAR_{I_i}\text{ (from $\XISTAR(S_i)=\XQSTAR(S_i)$)}.
\end{align*}
Observe that $\LSTAR[i-1] + \LISTAR_{I_i}$ is a lower bound for any vector of interventions on interval $J_i$ which proves the optimality of $\XBF^*$ for \PROBLEMA{i}.

Finally, equality $\LSTAR[i]=\LSTARAUGMN[i-1]+\LISTAR_{I_i}$ straightly derives from equality $\LVECTBKW[x^*]{i}=\LSTAR[i-1] + \LISTAR_{I_i}$, optimality of $\bar{\XBF}$ for \PROBLEMA{i} and equation \EQUREF{PN14_37:augmntEqualNormal}.
\end{myproof}

\bigskip
\begin{myproof}[Proof of \TEOREF{PN:algoBKWcorrectness}]

The algorithm starts with some initialization and computes vectors
$\SHBF$, $\lpBF$, $\lmBF$, and $\SBF$ representing, respectively, virtual stock, surplus loss, stockout loss, and stock level induced by the null vector of interventions $\XBF=\OBF$ on all epochs $T$ (lines \LINEREF{TH3_1}--\LINEREF{TH3_8}).
The computation is clearly done in $O(m)$ time.

\medskip
In the second part, three fundamental steps are performed.
\begin {enumerate}

\item Lines \LINEREF{TH3_9} -- \LINEREF{TH3_11}.
A \ONEINTPROBLEM\ problem is solved in $O(m-e_\NUMV)$ time by procedure \CiteAlgoVehicleIntervention\ on the interval $I_\NUMV$ with initial stock level $\INISTOCKEPOCHBYVECT[\OBF]{\NUMV}$ and producing corresponding minimum modulus optimal interventions \XISTAR, \XQSTAR\ and optimal values \LISTAR, \LQSTAR\ for the uncapacitated and capacitated problem, respectively; the stock level at epoch $m$ is $s^{(m)}(\XQSTAR)=s^{(m)}(\OBF)$.

\item Line \LINEREF{TH3_12}.
The augmented problem $\PROBLEMAAUG{\NUMV-1}$ is solved with $\delta=\XISTAR-\XQSTAR$ in $O(e_\NUMV)$ time by algorithm \CiteAlgoGlobalBackwardRec\ producing an optimal vector of interventions \XBF\
and the corresponding optimal value $L$.

\item Lines \LINEREF{TH3_13} -- \LINEREF{TH3_14}.
The solutions of the two problems are chained in $O(1)$ time to form the optimal vector of interventions $\XBF^*$ for \PROBLEMA{\NUMV} such that
$s^{(m)}(\XBF^*)=s^{(m)}(\OBF)$.
and the corresponding optimal value $L^*$.
\end{enumerate}

Correctness of the first step is guaranteed by  (\TEOREF{TH2:alg:Vehicle-Intervention}); correctness of the third step is guaranteed by  (\PROREF{PN14:MERGE}). Correctness of the second step is proved in the following.

\medskip
We proceed by induction on the number $i$ of stages to show that problem \PROBLEMAAUG{i} is solved by algorithm \CiteAlgoGlobalBackwardRec\ which produces an optimal vector of interventions $\XBF^*$\ in $O(e_{i+1})$ time.

\medskip\noindent \textbf{Base cases}

For $i=0$, we do not have any decision epochs and the procedure consistently returns in $O(1)$ time the empty vector as optimal vector of interventions $\XBF^*$ and $|\delta|$ as optimal value $L^*$.

For $i=1$, the augmented problem \PROBLEMAAUG{1} reduces to a \ONEINTPROBLEM\ problem on augmented interval $\AUGMINTI[1]$ which is consistently initialized at lines \LINEREF{TH3R_2} -- \LINEREF{TH3R_5} in $O(1)$ time and solved by algorithm \CiteAlgoVehicleIntervention\ at line \LINEREF{TH3R_6} in $O(e_{i+1})$, which returns the corresponding optimal interventions \XQSTAR, \XISTAR, and optimal values \LQSTAR\ and \LISTAR.
The next recursive call for \PROBLEMAAUG[\XISTAR-\XQSTAR]{0}, returns an empty vector of interventions and a loss $|\XISTAR-\XQSTAR|$.
Thus, the procedure consistently returns  $\XBF^*=[\XISTAR]$ as one-intervention optimal vector of interventions and $L^*=\LQSTAR=\LISTAR+|\XISTAR-\XQSTAR|$ as total loss (latter equality guaranteed by equation \EQUREF{PN7_16:PALGO.4.eq3} of \PROREF{PN7:XQfromXI}).

\medskip\noindent \textbf{Inductive step.}

Let $i>1$ and assume that
algorithm \CiteAlgoGlobalBackwardRec\ solves problem \PROBLEMAAUG{i-1} in $O(e_{i})$ time. We show that algorithm \CiteAlgoGlobalBackwardRec\ also solves problem \PROBLEMAAUG{i} in $O(e_{i+1})$.

Similarly to the case of $i=1$, also for $i>1$,
data for the \ONEINTPROBLEM\ problem on augmented interval $\AUGMINTI$ are computed in $O(1)$ and corresponding problem is solved in $O(e_{i+1}-e_i)$ by algorithm \CiteAlgoVehicleIntervention, which returns the corresponding optimal interventions \XQSTAR, \XISTAR, and optimal values \LQSTAR\ and \LISTAR.

From the inductive hypothesis, augmented problem \PROBLEMAAUG[\XISTAR-\XQSTAR]{i-1} is solved by algorithm \CiteAlgoGlobalBackwardRec\ in $O(e_{i})$ producing an optimal vector of interventions \XBF\ and its optimal value $L$.

According to \PROREF{PN14:MERGE}, at last step the optimal vector of interventions $\XBF^*$ for the augmented problem \PROBLEMAAUG{i} is given by $[\XBF,\XQSTAR]$ and its optimal value is given by $L^*=\LISTAR+L$. The overall time complexity is $O(e_{i+1})$.
\end{myproof}

\subsection{Table of notation}

\footnotesize
\begin{longtable}{p{3.8cm}p{12.5cm}}
\hline
        \textbf{Notation} & \textbf{Meaning} \\ \hline \endhead
        \textbf{Problem parameters} \\ \hline
        $C$ & Capacity of the station. \\
        $T=\{1,\ldots,m\}$ & Set of epochs. \\
        $\DBF=[d^{(1)},\ldots,d^{(m)}]$ & Vector of net flow in the time horizon. \\
        $V=\{v_1,\ldots,v_\NUMV\}$ & Set of the vehicles. \\
        $H=\{e_1,\ldots,e_\NUMV\}\subseteq T$ & Set of epochs in which the station is visited by the vehicles. \\
        $\QBF=[Q_1,\ldots,Q_\NUMV]$ & Vector of the capacity of vehicles. \\
        $\qBF=[q_1,\ldots,q_\NUMV]$ & Vector of the bike load of the vehicles. \\
        \hline
        \textbf{Variables} \\ \hline
        $\XBF=[x_1,\ldots,x_\NUMV]$ & Vector of the interventions. \\
        $\SBF=[s^{(1)},\ldots,s^{(m)}]$ & Vector of stock levels in the time horizon. \\
        $\SHBF = [\hat{s}^{(1)},\ldots,\hat{s}^{(m)}]$ & Vector of virtual stock in the time horizon. \\
        $\lpBF = [l^{(1)+},\ldots,l^{(m)+}]$ & Vector of surplus losses. \\
        $\lmBF = [l^{(1)-},\ldots,l^{(m)-}]$ & Vector of stockout losses. \\
        $L$ & Total amount of lost requests. \\
        \hline
        \textbf{Interval optimization} \\
        \hline
        $\STARTEPOCH\in T$ & Starting epoch for interval optimization. \\
        $\ENDEPOCH\in T$ & Ending epoch for interval optimization. \\
        $\SMINT{h}(x)$ & Minimum stock level in epoch interval $[\STARTEPOCH, h]\subseteq T$. \\
        $\SMAXT{h}(x)$ & Maximum stock level in epoch interval $[\STARTEPOCH, h]\subseteq T$. \\
        $h^{(+)}(x)$ & First epoch in which a surplus loss is present, if any, $+\infty$ otherwise. \\
        $h^{(-)}(x)$ & First epoch in which a stockout loss is present, if any, $+\infty$ otherwise. \\
        $\XISSTAR=\XIRSTAR$ & Range of optimal interventions for the uncapacitated case. \\
        $\LISTAR$, $\XISTAR$ & Optimal loss and minimum modulus optimal intervention for the uncapacitated case. \\
        $\XQSSTAR=\XQRSTAR$, $\LQSTAR$, $\XQSTAR$ & Same as above for the capacitated problem. \\
        \hline
        \multicolumn{2}{p{16cm}}{\textbf{Time horizon optimization}} \\ \hline
        $I_i$ & $I_i=[e_i,e_i+1,\ldots,e_{i+1}-1] \text{ for } i=1,\ldots,\NUMV-1$, $I_\NUMV=[e_\NUMV,\ldots,m]$. \\
        $J_i$ & $J_i=[e_1,e_2,\ldots,e_{i+1}-1] \text{ for } i=1,\ldots,\NUMV-1$, $J_\NUMV=[e_1,\ldots,m]$. \\
        $\AUGMINTI$, $\AUGMINTJ$ & Augmented intervals $\AUGMINTI=I_i\cup\{e_{i+1}\}$, $\AUGMINTJ=J_i\cup\{e_{i+1}\}$. \\
        $\PROBLEMA{i}$ & 
        $\PROBLEMA{i}\DEFINE\PIOTTBKW{i}$, optimization problem \EQUREF{eq07:base:L} on interval $J_i$. \\
        $\PROBLEMAAUG{i}$ & $\PROBLEMAAUG{i}\DEFINE\PIOTTBKWAUGMN{i}$, the augmentation of problem \PROBLEMA{i} to \AUGMINTJ\ with extra loss $\delta$.\\ 
        \LSTAR, \LSTARAUGMN & Optimal values of \PROBLEMA{i} and \PROBLEMAAUG{i}. \\
        $\XVECTrow[x]{1}{i}$ & 
        Subvector of first $i$-th components of the vector of interventions \XBF. \\
        $\LVECTBKW{i}$ & Loss induced by interventions $\XVECTrow{1}{i}$ in interval in interval $J_i$.\\ 
        $\LVECTBKWAUGMN{i}$ & Loss induced by interventions $\XVECTrow{1}{i}$ in interval \AUGMINTJ\ for the augmented problem. \\ 
        $\LIXS[x]{i}{s}$ & Loss induced by intervention $x$ in interval $I_i$ depending on initial stock level $s$ at epoch $e_i-1$. \\
        \hline
        \textbf{Rich notation} \\ \hline
        \multicolumn{2}{p{16cm}}{Throughout the manuscript, the dependency of different quantities on different variables is highlighted. Example of this highlighting through notation are reported in this section of the table.}\\ \hline
        $s^{(h)}(x)$ & Stock level at epoch $h$ resulting from intervention $x$ (for \ONEINTPROBLEM\ problem). \\
        $\SBF(\XBF)$ & Vector of stock levels resulting from the vector of interventions $\XBF$. \\
        $s^{(h)}(\XBF)$ & Stock level at epoch $h$ induced by interventions in vector $\XBF$ up to epoch $h$. \\
        $L_{[\STARTEPOCH,h]}(x\ |\ s)$ & Loss induced on interval $[\STARTEPOCH,h]$ by intervention $x$ at epoch $\STARTEPOCH$ given an initial stock level $s$. \\
        $\XISSTAR_{I_k}(s)=[\XIDSTAR_{I_k}(s),\XIUSTAR_{I_k}(s)]$, $\XISTAR_{I_k}(s)$, $\LISTAR_{I_k}(s)$ & Set of optimal interventions, minimum modulus optimal intervention and optimal loss for the \ONEINTPROBLEM\ uncapacitated problem on interval $I_k$, with their dependency on initial stock level $s$ being highlighted. \\
        \hline
        \\
    \caption{Table of notation}
    \label{tab:my_label}
\end{longtable}

\normalsize

\newpage

\bibliographystyle{abbrvnat}
\bibliography{biblio}
\end{document}